\documentclass[11pt]{article}
\usepackage{amssymb}
\usepackage{amsfonts}
\def\build#1_#2^#3{\mathrel{\mathop{\kern 0pt#1}\limits_{#2}^{#3}}}

\newdimen\emm
\def\pmb#1{\emm=0.03em\leavevmode\setbox0=\hbox{#1}
\kern0.901\emm\raise0.434\emm\copy0\kern-\wd0
\kern-0.678\emm\raise0.975\emm\copy0\kern-\wd0
\kern-0.846\emm\raise0.782\emm\copy0\kern-\wd0
\kern-0.377\emm\raise-0.000\emm\copy0\kern-\wd0
\kern0.377\emm\raise-0.782\emm\copy0\kern-\wd0
\kern0.846\emm\raise-0.975\emm\copy0\kern-\wd0
\kern0.678\emm\raise-0.434\emm\copy0\kern-\wd0
\kern\wd0\kern-0.901\emm}

\newcommand {\irm}{{I}_{\rho , \mu}}
\newcommand {\irmp}{{I}_{\rho , \mu}^+}

\newcommand {\ikl}{{I}_{\kappa, \lambda}}
\newcommand {\iol}{{I}_{\omega, \lambda}}

\newcommand {\ikle}{{I}_{\kappa, \lambda}^*}
\newcommand {\iklp}{{I}_{\kappa, \lambda}^+}

\newcommand {\prmu}{P_{\rho}(\mu)}
\newcommand {\pol}{P_{\omega}(\lambda)}

\newcommand {\pkl}{P_{\kappa}(\lambda)}
\newcommand {\ptd}{P_{\theta}(\delta)}
\newcommand {\pait}{P_{|a\cap\theta|}(a)}
\newcommand {\dt}{{\lbrack\delta\rbrack^{<\theta}}}
\newcommand {\lk}{{\lbrack\lambda\rbrack^{<\kappa}}}
\newcommand {\jpia}{J^+\cap P(A)}

\newcommand {\nsdt}{NS_{\kappa,\lambda}^{\lbrack\delta\rbrack^{<\theta}}}
\newcommand {\nsdl}{NS_{\kappa,\lambda}^{\lbrack\lambda\rbrack^{<\theta}}}
\newcommand{\dom}{\mathfrak{d}}

\newcommand{\deklmu}{\delta_{\kappa,\lambda}^\mu}

\newcommand{\dklmu}{\dom_{\kappa,\lambda}^\mu}

\newcommand{\emupkl}{{^\mu}\!{\pkl}}
\newcommand{\emupkk}{{^\mu}\!{P_\kappa(\kappa)}}

\newcommand{\emupka}{{^\mu}\!{P_\kappa(\alpha)}}

\newcommand {\eop}{\hfill$\square$\\} 
\newcommand {\boite}{\square}

\begin{document}


\thispagestyle{empty}
\begin{center}
       {\Large{\sc{COFINALITY OF NORMAL IDEALS ON $\pkl$}}\quad I}
\end{center}
\vspace*{2cm}

\begin{center}
{\sc Pierre Matet}\\
\vspace*{2ex}
  Universit\'e de Caen-CNRS\\
Math\'ematiques\\
BP 5186\\
14032 Caen Cedex, France\\
matet@math.unicaen.fr\\
\vspace*{0.5cm}
{\sc C\'edric P\'ean}
\footnote{Some of the material in this paper originally
appeared as part of the author's doctoral
\mbox{dissertation} completed at the Universit\'e de Caen, 1998.}\\
\vspace*{2ex}
  Universit\'e de Caen-CNRS\\
Math\'ematiques\\
BP 5186\\\
14032 Caen Cedex, France\\
\vspace*{0.5cm}
{\sc Saharon Shelah}
\footnote{Partially supported by the Israel Science
   Foundation. Publication 713.
\bigskip

2000 Mathematics Subject Classification~: 03E05, 03E55, 03E35.

}\\
\vspace*{2ex}
  Institute of Mathematics\\
The Hebrew University of Jerusalem\\
91904 Jerusalem, Israel\\
\vspace*{2ex}
Department of Mathematics\\
Rutgers University\\
New Brunswick, NJ 08854,  USA\\
shelah@math.huji.ac.il
\vskip3cm
\end{center}

\baselineskip=22pt

\newpage
\thispagestyle{empty}
\begin{abstract}
\noindent Given an ordinal $\delta\le\lambda$ and a cardinal
$\theta\le\kappa$, an ideal $J$ on $\pkl$ is said to be
$\dt$-normal if given $B_e\in J$ for  $e\in P_\theta(\delta)$,  the 
set of all $a\in\pkl$ such
that $a\in B_e$ for some $e\in P_{|a\cap\theta|}(a\cap\delta)
$ lies in $J$.
We give necessary and sufficient
conditions for the existence of such ideals and describe the least 
one, denoted by $\nsdt$.
We compute the cofinality of $\nsdt$.
\end{abstract}


\medskip
{\bf\Large 0. Introduction}
\smallskip

Given a regular uncountable cardinal $\kappa$ and a cardinal 
$\lambda\ge\kappa$, an ideal on $\pkl$
is said to be normal if it is closed under diagonal unions of $\lambda$ many
of its members. Building on work of Jech [J] and Menas [Me], Carr [C] described
the least such ideal, usually denoted by $NS_{\kappa,\lambda}$. 
Numerous variations
on the original notion of normality have been considered  over the 
years. We are
interested in two of these variants. First, there is a notion called  `strong
normality' which has been rather extensively studied (see e.g. [CP], [F],
[M1], [CLP]). The definition involves diagonal unions of length
$\lambda^{<\kappa}$.
[CLP] gives necessary and sufficient conditions for the existence of strongly
normal  ideals and describes the least such ideal when there is one. As the
terminology implies, every strongly
normal ideal is normal. The other notion is that of $\delta$-normality for an
ordinal $\delta\le\lambda$.  An ideal on $\pkl$ is called $\delta$-normal if it
is closed under diagonal unions of length $\delta$.  Thus 
$\lambda$-normality is the
same as normality.
$\delta$-normality has been  studied by Abe [A] who gave a 
description of the smallest
$\delta$-normal ideal on $\pkl$.

\smallskip
We introduce a more general concept, that of \ $\dt$-normality, where
$\delta$ is, as above, an ordinal with
$\delta\le\lambda$, and $\theta$ is a cardinal with 
$\theta\le\kappa$. The definition
is similar to that of strong normality,  with this difference that our diagonal
unions are indexed by $\dt$. So $\lk$-normality is  identical with 
strong normality,
whereas
$\lbrack\delta\rbrack^{<2}$-normality is the same as
$\delta$-normality.

\smallskip
We give necessary and sufficient conditions for the existence of 
$\dt$-normal ideals
on $\pkl$ and
describe the least such ideal, which we denote by $\nsdt$.

$[\lambda]^{<\theta}$-normality \ (for  $\theta$ a regular infinite 
cardinal less than $\kappa$)
has been independently considered by D\v zamonja [D]. In particular, 
Claim 2.9 and Corollary
2.13 of [D] provide alternative descriptions of \ $\nsdl$.

\smallskip
Given an ideal $J$, its cofinality $cof(J)$ is its least number of 
generators, i.e. the least
size of any subcollection $X$ of $J$ such that every member of the 
ideal is included in some
element of $X$. We determine the cofinality of $\nsdt$. Its computation
involves a multidimensional version of  the dominating number
$\dom_\kappa$, which is no surprise, as Landver
$\lbrack$Lemma 1.16 in  [L])
proved that the cofinality of the minimal normal ideal on $\kappa$ is 
$\dom_\kappa$.

Part of the paper is concerned with the problem of comparing the 
various ideals that
are considered.  Given two pairs
$(\delta,\theta)$ and $(\delta^\prime,\theta^\prime)$, we investigate 
whether $\nsdt$ and
$NS_{\kappa,\lambda}^{\lbrack\delta^\prime\rbrack^{<\theta^\prime}}$ 
are equal, and,
more generally, whether one of the two ideals is a restriction of the other
(there is more about this in [MP\'eS]). It is
for instance shown that
$\nsdt=NS_{\kappa,\lambda}^{\lbrack|\delta|\rbrack^{<\theta}}|A$
for some $A$.

\smallskip
Section 1 collects basic definitions and facts concerning ideals on $\pkl$.
This is \mbox{standard}
  material except
for Proposition 1.4. In Section 2 we introduce the property
of \break $\dt$-normality and state necessary
and sufficient conditions for the existence of a \break\mbox{$\dt$-normal}
ideal on $\pkl$.
The discussion is very much like the one regarding the existence of a
strongly normal ideal, and
arguments are routine. We briefly consider various weaker properties
(compare Proposition 2.3
((iii) and (iv)), Proposition 2.6 (ii) and Corollary 2.8 (ii)
  with Proposition 2.5 (ii))
and characterize the ideals that satisfy them. In Sections 3 and 4 we show
that we could without loss of generality
assume that $\theta$ is an infinite cardinal and $\delta$ is either a cardinal
less than $\kappa$,  or a multiple of $\kappa$.
We describe the smallest $\dt$-normal ideal on $\pkl$, denoted by
$\nsdt$.  Section 5 isconcerned with the case that $\theta$ is a limit
cardinal. It is proved that if $\delta\geq\kappa$
and $\theta$ is a singular strong limit cardinal, then every 
$[\delta]^{<\theta}$-normal ideal on
$P_\kappa(\lambda)$ is $[\delta]^{<\theta^+}$-normal. Sections 6 and 
7 deal with the question of
the existence of an ordered pair $(\delta', \theta') \not= 
(\delta,\theta)$ such that
$\delta'\leq\delta$, $\theta'\leq\theta$ and 
$NS_{\kappa,\lambda}^{[\delta]^{<\theta}} =
NS_{\kappa,\lambda}^{[\delta']^{<\theta'}} | A$ for some $A$.

\smallskip
In Section 8 we introduce a three-cardinal version, denoted by $\dklmu$,
of the dominating number
$\dom_\kappa$. There are many identities involving the $\dklmu$'s, 
and we present some of them. The
cofinality of $\nsdt$ is computed in Section 9.

\vfill\eject
{\Large\bf 1. Ideals}\\

{\bf Definition.} \ {\it Given a set A and a cardinal $\tau$, we put 
$P_\tau(A)=\lbrack
A\rbrack^{<\tau}=\{a\subseteq A :|a|<\tau\}$.}

\bigskip
{\bf Throughout the section $\rho$  will denote an infinite cardinal 
and $\mu$ a
cardinal with  $\mu\ge\rho$}.

\bigskip
The section presents some basic material concerning ideals on $\prmu$. Let
us start by recalling some definitions.\\

{\bf Definition.} \ {\it We set $\widehat{a}=\{b\in\prmu:a\subseteq 
b\}$ for every
$a\in\prmu$.}

\bigskip
{\bf Definition.} \  {\it $\irm$ is the collection of all 
$B\subseteq\prmu$ such
that $B\cap\widehat{a}=\emptyset$ for some $a\in\prmu$.}

\bigskip
{\bf Definition.} \ {\it By an {\it ideal on} $\prmu$, we mean a 
collection $K$ of subsets
of $\prmu$ such that (0) $P (B)\subseteq K$ for all $B\in K$, (1) 
$\cup Y\in K$ for
all $Y\subseteq K$ with
$0<|Y|<cf(\rho)$, (2) $\irm\subseteq K$, and (3) $\prmu\not\in K$.}

\bigskip
{\bf Definition.} \ {\it Two ideals $I,J$ on $\prmu$ {\it cohere} if 
$I\cup J\subseteq K$
for some ideal $K$ on $\prmu$.}

\medskip
The following is easily verified.\\

\bigskip
\hskip-0,2cm\pmb{\sc Proposition 1.1.} \  {\it $\irm$ is an ideal on $\prmu$.}

\bigskip
{\bf Definition.} \ {\it Let $K$ be an ideal on $\prmu$.

We set $K^+=P(\prmu)-K$ and $K^*=\{B\subseteq\prmu:\prmu-B\in K\}$.

$non(K)$ is the least cardinality of any $A\subseteq \prmu$ with $A\in K^+$.

$cof(K)$ is the least cardinality of any $S\subseteq K$ with
$K=\displaystyle\bigcup_{B\in S}P(B)$.}
\bigskip

The following is well-known. \\

\bigskip
\hskip-0,2cm\pmb{\sc Proposition 1.2.} \  {\it Let $K$ be an ideal on $\prmu$.
Then $non(K)\le cof(K)$}.\\
\vspace*{2ex}{\it Proof. }Let $S\subseteq K$ be such that
$K=\displaystyle\bigcup_{B\in S}P(B)$. Pick $a_B\in\prmu-B$  for  $B\in S$.
Then $\{a_B:B\in S\}\in K^+$.\eop

{\bf Definition.} \ {\it We put $u(\rho,\mu)=non(\irm)$}.\\

\hskip-0,2cm\pmb{\sc Proposition 1.3.} \
\begin{enumerate}
\item[(0)] \vskip-0,3cm $\mu\le u(\rho,\mu)$.

\item[(1)] \vskip-0,3cm $cf(\rho)\le cf(u(\rho,\mu))$.
\end{enumerate}

\vspace*{0ex}{\it Proof. }
\vskip-0,3cm
\begin{enumerate}
\vskip-0,5cm
\item[(0)] : Given $A\in\irmp$, we have $\mu=\cup
A$ and  therefore $\mu\le\rho\cdot|A|$.
This proves the desired inequality if  $\mu>\rho$.
Given $B\subseteq{P}_\rho(\rho)$ with $|B|<\rho$, pick 
$\alpha_b\in\rho-b$ for  $b\in B$.
Then $\{\alpha_d:d\in B\}\not\subseteq b$ for all $b\in B$, and 
consequently $B\in
I_{\rho,\rho}$.  Hence $u(\rho,\rho)\ge\rho$.

\item[(1)] : Use the fact that $\prmu$ is closed under unions of less 
than $cf(\rho)$ many
of its\break  members.~\eop
\end{enumerate}
The following result will be used in Section 8.\\

\hskip-0,2cm\pmb{\sc Proposition 1.4.} \  {\it Let $K$ be an ideal on 
$\prmu$. Further let
$A\in K^+$, and set
$\chi=min\{|A\cap C |:C\in K^*\}$. Assume that $cof(K)\le \chi$. Then 
$\chi$ is the largest
cardinal $\tau$ such that there  exists a partition of $A$
into $\tau$ sets in $K^+$.} \\

\vskip-0,5cm
\vspace*{0ex}{\it Proof. }Pick $C\in K^*$ with $|A\cap C|=\chi$, and
set $D=A\cap C$.

\bigskip
Let us first suppose that there exists $g:A\longrightarrow\chi^+$ 
such that $g^{-1}(\{\alpha\})\in K^+$
  for all $\alpha\in \chi^+$. Then $\{C\cap 
g^{-1}(\{\alpha\}):\alpha\in\chi^+\}$ is a partition of $D$
  into $\chi^+$ pieces in $K^+$, which contradicts the fact that $|D|=\chi$.

Let us now show that there exists a partition of $A$ into $\chi$ sets 
in $K^+$. Select a bijection
$j:\chi\times\chi\longrightarrow\chi$, and let $B_\beta$ for 
$\beta<\chi$ be such
that $K=\displaystyle\bigcup_{\beta<\chi}P(B_\beta)$. Define $H_\xi$
for $\xi<\chi$ so that $H_{j(\alpha,\beta)}=
B_\beta$ for every $(\alpha,\beta)\in\chi\times\chi$. Notice that 
given $\beta,\eta<\chi$, there is
$\xi\ge\eta$ with $H_\xi=B_\beta$. Now construct $a_\eta^\xi$ for 
$\eta\le\xi<\chi$ so that

\begin{enumerate}
\item[(0)] \vskip-0,2cm
$\{a_\gamma^{\xi^\prime}:\gamma\le\xi^\prime\}\cap\{a_\eta^\xi:\eta\le 
\xi\}=\emptyset$ for
$\xi^\prime<\xi<\chi$.

\item[(1)] \vskip-0,2cm $a_\eta^\xi\not=a_{\eta^\prime}^\xi$ for 
$\eta^\prime<\eta\le\xi<\chi$.

\item[(2)] \vskip-0,2cm  $a_\eta^\xi\in D-(H_\xi\cup B_\eta)$ for 
$\eta\le\xi<\chi$.
\end{enumerate}

Set $A_\eta=\{a_\eta^\xi:\eta\le\xi<\chi\}$ for $\eta<\chi$. Then the 
following hold:
\begin{enumerate}
\item[(i)] \vskip-0,2cm $|A_\eta|=\chi$.

\item[(ii)] \vskip-0,2cm $A_\eta\in K^+\cap P(A)$.

\item[(iii)] \vskip-0,2cm $A_\eta\cap B_\eta=\emptyset$.

\item[(iv)] \vskip-0,2cm $A_\eta\cap A_{\eta^\prime}=\emptyset$ for 
$\eta^\prime<\eta$.\eop
\end{enumerate}

\vskip-0,3cm
\pmb{\sc Corollary 1.5.} \  {\it There exist $A_e\in\irmp\cap 
P(\widehat{e})$ for
$e\in\prmu$ such that  {\rm(a)}  $|A_e|=\mu^{ <\rho}$ for every 
$e\in\prmu$, and {\rm(b)}
$A_e\cap A_{e^\prime}=\emptyset$ for all $e,e^\prime\in\prmu$
  with $e\not=e^\prime$.} \\
\vskip-0,3cm
\vspace*{0ex}{\it Proof. }By the proof of Proposition 1.4.\eop

{\bf Definition.} \ {\it Given an ideal $K$ on $\prmu$, we put 
$K|A=\{B\subseteq\prmu:B\cap A\in
K\}$ for every $A\in K^+$.}

\bigskip
\hskip-0,2cm\pmb{\sc Proposition 1.6.} \  {\it Let $K$ be an ideal on 
$\prmu$, and $A\in
K^+$. Then
$K|A$ is an ideal on
$\prmu$.  Moreover, $K\subseteq K|A$ and $cof(K|A)\le cof(K)$.}\\
\vskip-0,3cm
\vspace*{0ex}{\it Proof. }Use the fact that for every
$B\subseteq\prmu$, $B\in K|A$ if and only if
$B\subseteq E\cup(
\prmu-A)$ for some $E\in K$.\eop

We will make use of the following observation.\\

\hskip-0,2cm\pmb{\sc Proposition 1.7.} \  {\it Let $I, J, K$ be three 
ideals on $\prmu$ such
that
$I\subseteq J\subseteq K$.
  Assume that there exists $A\in I^+$ such
that $K=I|A$. Then $J|A=I|A$.} \\
\vspace*{2ex}{\it Proof. }Notice that $A\in J^+$ since $A\in K^*$.
For each $B\subseteq\prmu$, we have
\vskip-0,3cm
\centerline{$B\in I|A\Rightarrow B\in J|A$}
\vskip-0,3cm
and
\vskip-0,3cm
\centerline{$B\not\in I|A\Rightarrow B\cap A\not\in I|A\Rightarrow B
\cap A\not\in K\Rightarrow B\not\in J|A$.}\eop


{\Large\bf 2. $\lbrack\delta\rbrack^{<\theta}$-normality}\\

{\bf Throughout the remainder of the paper $\kappa$  denotes a 
regular infinite cardinal,
$\lambda$ a cardinal with \ $\lambda \geq \kappa$, \ $\theta$ a cardinal with
$2\leq\theta\leq\kappa$, and $\delta$ an ordinal with 
$1\leq\delta\leq\lambda$}.

\bigskip
{\bf We set $\overline{\theta}=\theta$ if $\theta<\kappa$, or 
$\theta=\kappa$ and
$\kappa$ is a limit cardinal, and
$\overline{\theta}=\nu$ if $\theta=\kappa=\nu^+$}.

\bigskip
{\bf Throughout the remainder of the paper $J$  denotes a fixed ideal on
$\pkl$}.

\bigskip
In this section we introduce the notion of 
$\lbrack\delta\rbrack^{<\theta}$-normal ideal on $\pkl$
and describe necessary and sufficient conditions for the existence of 
such ideals. Let us start with a
few definitions.\\

{\bf Definition.} \ {\it Given $X_e\subseteq\pkl$ for $e\in\ptd$, we let
\vskip-0,5cm
$$\build{\Delta}_{e\in\ptd}^{} 
X_e=\bigcap_{e\in\ptd}(X_e\cup\{a\in\pkl:e\not\in
\pait\})$$
\vskip-0,5cm
and
\vskip-0,5cm
$$\build{\nabla}_{e\in\ptd}^{}X_e=\displaystyle\bigcup_{e\in\ptd}(X_e
\cap\{a\in\pkl:e\in
\pait\}).$$}

{\bf Definition.} \ {\it We define $\nabla^{\dt}J\subseteq P(\pkl)$ 
by : $B\in\nabla^{\dt}J$
if and only if there are $B_e\in J$ for $e\in\ptd$ such that
$$B\subseteq\{a\in\pkl:a\cap\theta=\emptyset\}
\cup(\build{\nabla}_{e\in\ptd}^{} B_e).$$}

\hskip-0,2cm\pmb{\sc Lemma 2.1.} \
\vskip-0,3cm
\begin{enumerate}
\item[(0)] \vskip-0,3cm $J\subseteq\nabla^{\dt}J$.
\item[(1)] \vskip-0,3cm $\cup Y\in\nabla^{\dt}J$ {\it for all} $Y\in 
P_\kappa(J)-\{\emptyset\}$.
\item[(2)] \vskip-0,3cm {\it Assume that $\delta^\prime$ is an ordinal with
$\delta\leq\delta'\leq\lambda$, \
$\theta^\prime$ is a  cardinal with  $\theta\leq\theta'\leq\kappa$, 
and $J^\prime$ is an
ideal on $\pkl$ with $J\subseteq J^\prime$. Then
$\nabla^{\dt}J\subseteq\nabla^{\lbrack\delta^\prime
\rbrack^{<\theta^\prime}}J^\prime$.}
\end{enumerate}
 
\vspace*{0ex}{\it Proof. }
\vskip-0,3cm
\begin{enumerate}
\item[(0)] :  It suffices to observe that
$B\subseteq\{a\in\pkl:a\cap\theta=\emptyset\}\cup 
(\build{\nabla}_{e\in\ptd}^{}B)$ for every $B\in J$.
\item[(1)]:   Use the fact that if  $X_e^\alpha\subseteq\pkl$ for 
$e\in\ptd$ and $\alpha<\rho$,
and $\rho$ is a cardinal with $\rho > 0$, then 
$\displaystyle\bigcup_{\alpha<\rho}(\build{\nabla}_{e\in\ptd}^{}X_e^\alpha)=
\build{\nabla}_{e\in\ptd}^{} (\displaystyle\bigcup_{\alpha<\rho}X_e^\alpha)$.
\item[(2)]:   Use (0),(1) and the fact that 
$\build{\nabla}_{e\in\ptd}^{} B_e\in
\nabla^{\lbrack\delta^\prime\rbrack^{<\theta^\prime}}J^\prime$ 
whenever $B_e\in J$
\item[ ]  for $e\in\ptd$.\eop
\end{enumerate}

\vskip-0,3cm
\hskip-0,2cm\pmb{\sc Proposition 2.2.}
\vskip-0,3cm
\begin{enumerate}
\item[\rm (0)] \vskip-0,2cm
$\nabla^{\dt}J=\nabla^{\lbrack\delta\rbrack^{<\overline{\theta}}}J$.
\item[\rm (1)] \vskip-0,2cm \vskip-0,3cm If $|\ptd|<\kappa$, {\it 
then} $J=\nabla^{\dt}J$.
\end{enumerate}

{\it Proof. }
\vskip-0,3cm
\begin{enumerate}
\vskip-0,3cm
\item[(0)] : Assume that $\theta=\kappa=\nu^+$.
Then  clearly, $P(\widehat{\nu})\cap
\nabla^{\lbrack\delta\rbrack^{<\kappa}}J=P(\widehat{\nu})\cap\nabla^{\lbrack
\delta\rbrack^{<\nu}}J$. Hence $\nabla^{\lbrack\delta\rbrack^{<\kappa}}J=
\nabla^{\lbrack\delta\rbrack^{<\nu}}
J$ by Lemma 2.1 $((0)$ and $(1))$.
\item[(1)]  : Use Lemma 2.1 $(0)$.\eop
\end{enumerate}
 
{\bf Definition.} \ {\it Given $A\subseteq\pkl$, 
$f:A\longrightarrow\ptd$ is {\it
$\ptd$-regressive} if
$f(a)\in\pait$ for  all $a\in A$ with $a\cap\theta\not=\emptyset$.}\\

\eject
\hskip-0,2cm\pmb{\sc Proposition 2.3.} \ {\it The following are equivalent :}
\vskip-0,3cm
\begin{enumerate}
\item[(i)] \vskip-0,3cm $\pkl\not\in\nabla^{\dt}J$.

\item[(ii)] \vskip-0,2cm $\nabla^{\dt}J$ {\it is an ideal on}  $\pkl$.

\item[(iii)]\vskip-0,2cm  $\build{\Delta}_{e\in\ptd}^{}C_e\in J^+$ 
{\it whenever} $C_e\in J^*$ {\it
for}  $e\in\ptd$.

\item[(iv)]\vskip-0,2cm $\build{\Delta}_{e\in\ptd}^{}C_e\in\iklp$ 
{\it whenever} $C_e\in J^*$ {\it
for}
$e\in\ptd$.

\item[(v)] \vskip-0,2cm {\it For every $\ptd$-regressive 
$f:\pkl\longrightarrow\ptd$, there is
$D\in J^+$ such  that $f$ is constant on $D$.}
\end{enumerate}

\vspace*{0ex}{\it Proof. }
\begin{enumerate}
\item[(i)]  $\rightarrow$ (ii)  : By Lemma 2.1 $((0)$ and $(1))$.
\item[(ii)] $\rightarrow$ (iii) : Use Lemma 2.1 $(0)$ and the fact that
\vskip-0,9cm
$$\build{\Delta}_{e\in\ptd}^{}C_e=\pkl- \build{\nabla}_{e\in\ptd}^{} 
(\pkl-C_e)$$

whenever $C_e\subseteq\pkl$ for  $e\in\ptd$.

\item[(iii)] $\rightarrow$ (iv) : Trivial.
 
\item[(iv)] $\rightarrow$ (v) : Use the fact that
\vskip-1,2cm
$$\build{\Delta}_{e\in\ptd}^{}(\pkl-f^{-1}(\{e\}))\in \ikl$$
for every $\ptd$-regressive $f:\pkl\longrightarrow\ptd$.
\item[(v)] $\rightarrow$ (i) : Assume that there are $B_e\in J$ for 
$e\in\ptd$ such that
\vskip-1cm
$$\{a\in\pkl:
a\cap\theta\not=\emptyset\}\subseteq\build{\nabla}_{e\in\ptd}^{}B_e.$$
Then there is a $\ptd$-regressive $f:\pkl\longrightarrow\ptd$ with 
the property that $a\in
B_{f(a)}$  for all $a\in\pkl$ with $a\cap\theta\not=\emptyset$. 
Clearly, $f^{-1}(\{e\})\in J$ for
every $e\in\ptd$.
\end{enumerate}\eop

{\bf Definition.} \ {\it $J$ is {\it $\dt$-normal} if $J=\nabla^\dt J$.}\\

\bigskip
\hskip-0,2cm\pmb{\sc Proposition 2.4.} \ {\it Let $\delta'$ be an 
ordinal with $1\leq
\delta'\leq\delta$, and $\theta'$ a cardinal with 
$2\leq\theta'\leq\theta$. Then every
$[\delta]^{<\theta}$- normal ideal on $P_\kappa(\lambda)$ is 
$[\delta']^{<\theta'}$ -
normal.}

\vspace*{2ex}{\it Proof. } By Lemma 2.1(2). \eop

\hskip-0,2cm\pmb{\sc Proposition 2.5.} \  {\it The following are equivalent:}
\begin{enumerate}
\item[(i)] $J$ {\it is} $\dt$-normal.

\item[(ii)] $\build{\Delta}_{e\in\ptd}^{} C_e\in J^*${\it  whenever 
$C_e\in J^*$ for
$e\in\ptd$.}

\item[(iii)] $\pkl\not\in\nabla^\dt(J|A)$ {\it for}  $A\in J^+$.

\item[(iv)] {\it Given $A\in J^+$ and a $\ptd$-regressive 
$f:A\longrightarrow\ptd$, there is
$D\in\jpia$ such  that $f$ is constant on $D$.}
\end{enumerate}

\vspace*{0ex}{\it Proof. }
\vskip-0,4cm
\begin{enumerate}
\vskip-1cm
\item[(i)]$\leftrightarrow$ (ii) : Use Lemma 2.1
$(0)$.

\item[(iii)]$\leftrightarrow$ (iv) : By Proposition 2.3 
((i)$\leftrightarrow$(v)) and Lemma
2.1 (0).

\item[(iii)]$\rightarrow$ (ii) : Use Proposition 2.3  ((i)$\rightarrow$(iii)).
\item[(ii)]$\rightarrow$ (iii) : Use Proposition 2.3 
((iv)$\rightarrow$(i)) and the fact that
\vskip-0,9cm
$$A~\cap\build{\Delta}_{e\in\ptd}^{} 
X_e=A~\cap\build{\Delta}_{e\in\ptd}^{}((\pkl-A)\cup
X_e)$$
whenever $A\subseteq\pkl$ and $X_e\subseteq\pkl$ for  $e\in\ptd$.\eop
\end{enumerate}
\vskip-0,7cm
Proposition 2.6 ((i)$\leftrightarrow$(iii))   shows that the
$\dt$-normality of $J$ can be seen as a  global property which 
corresponds to the local
property ``$\;\pkl\not\in\nabla^\dt J\;$". Let us next  briefly 
consider a weaker (see
Corollary 2.7) local property. The corresponding global property 
will be dealt with in
Corollary 2.8.\\

\vfill\eject
\hskip-0,2cm\pmb{\sc Proposition 2.6.} \ {\it Assume 
$\pkl\not\in\nabla^\dt\ikl$. Then the
following are equivalent:}
\begin{enumerate}
\item[(i)] {\it $J$ and $\nabla^\dt\ikl$ cohere.

\item[\rm (ii)] $\build{\Delta}_{e\in\ptd}^{} C_e\in J^+$ whenever 
$C_e\in\ikle$ for
$e\in\ptd$.

\item[\rm (iii)] Given $A\in J^*$ and a $\ptd$-regressive 
$f:A\longrightarrow\ptd$, there is
$D\in\iklp\cap
  P(A)$ such that $f$ is constant on $D$.}
\end{enumerate}

\vspace*{0ex}{\it Proof. }
\vskip-0,5cm
\begin{enumerate}
\vskip-0,5cm
\item[(i)] $\rightarrow$ (ii) : Straightforward.

\item[(ii)] $\rightarrow$ (iii) : Let $A\in J^*$ and 
$f:A\longrightarrow\ptd$ with the property
that $f^{-1}(\{e\})\in
\ikl$ for  $e\in\ptd$. Then $f(a)\not\in\pait$ for all $a\in
A \ \cap\build{\Delta}_{e\in\ptd}^{} (\pkl-f^{-1}(\{e\}))$.

\item[(iii)] $\rightarrow$ (i) : Assume that (iii) holds. Given 
$B_e\in\ikl$ for $e\in\ptd$,
define\\
$f: \build{\nabla}_{e\in\ptd}^{} B_e\longrightarrow\ptd$ so that for every
$a\in\build{\nabla}_{e\in\ptd}^{} B_e$, $f(a)\in\pait$ and $a\in 
B_{f(a)}$. Then $f$ is $\ptd$-regressive.
Moreover, $f^{-1}(\{e\})\in\ikl$ for every $e
\in\ptd$. It follows that $\build{\nabla}_{e\in\ptd}^{} B_e\not\in 
J^*$. Hence, setting
$K=\{B\cup E:B\in J$ and $E\in\nabla^{\dt}\ikl\}$,
we have that $K$ is an ideal on $\pkl$ with 
$J\cup(\nabla^{\dt}\ikl)\subseteq K$.
\eop
\end{enumerate}

\hskip-0,2cm\pmb{\sc Corollary 2.7.} \  {\it If 
$\pkl\not\in\nabla^{\dt}J$, then $J$ and
$\nabla^{\dt}\ikl$ cohere.}\\
\vspace*{0ex}{\it Proof. }By Lemma 2.1 $(2)$, Proposition 2.3
((i) $\rightarrow$ (iii)) and Proposition 2.6 ((ii)
$\rightarrow$ (i)).\eop

\hskip-0,2cm\pmb{\sc Corollary 2.8.} \  {\it Assume 
$\pkl\not\in\nabla^{\dt}\ikl$. Then
the following are equivalent~:}
\begin{enumerate}
\item[(i)] {\it $J|A$ and $\nabla^{\dt}\ikl$ cohere for every $A\in J^+$.}

\item[(ii)]\vskip-0,2cm {\it $\build{\Delta}_{e\in\ptd}^{} C_e\in 
J^*$ whenever $C_e\in\ikle$ for
$e\in\ptd$.}

\item[(iii)] {\it Given $A\in J^+$ and a $\ptd$-regressive 
$f:A\longrightarrow\ptd$,
there is $D\in\iklp\cap P(A)$ such that $f$ is constant on $D$.}

\item[(iv)] \vskip-0,2cm {$\nabla^{\dt}\ikl\subseteq J$.}
\end{enumerate}

We will now show that $\lbrack\delta\rbrack^{<2}$-normality is the 
same as $\delta$-normality
(which was studied by Abe in [A]). Let us first recall the following 
definitions.\\

{\bf Definition.-} \ {\it Given $X_\alpha\subseteq\pkl$ for 
$\alpha<\delta$, we set
\vskip-0,5cm
$$\build{\Delta}_{\alpha<\delta}^{} X_\alpha=
\bigcap_{\alpha<\delta}(X_\alpha\cup(\pkl-\widehat{\{\alpha\}}))$$
\vskip-0,3cm
and
\vskip-0,5cm
$$\build{\nabla}_{\alpha<\delta}^{} 
X_\alpha=\bigcup_{\alpha<\delta}(X_\alpha\cap
\widehat{\{\alpha\}}).$$}

{\bf Definition.-} \ {\it Given \ $K\subseteq P(P_\kappa (\lambda))$, 
we define $\nabla^\delta
K\subseteq P(\pkl)$ by :
$B\in\nabla^\delta K$  if and only if there are $B_\alpha\in K$ for
$\alpha<\delta$ such that
$B\subseteq(\pkl-\widehat{\{0\}})\cup\build{\nabla}_{\alpha<\delta}^{} 
B_\alpha$.}

\bigskip
{\bf Definition.-} \ {\it $J$ is {\it $\delta$-normal} if 
$J=\nabla^\delta J$.}\\

\hskip-0,2cm\pmb{\sc Proposition 2.9.} \  {\it $J$ is $\delta$-normal 
if and only if $J$ is
$\lbrack
\delta\rbrack^{<2}$-normal.}\\
\vspace*{0ex}{\it Proof. } The result easily follows from the
following two remarks:
\begin{enumerate}
\item[1)] \vskip-0,3cm Let $X_\alpha\subseteq\pkl$ for 
$\alpha<\delta$. Define $Y_e$ for $e\in
P_2(\delta)$ by : $Y_{\{\alpha\}}=X_\alpha$ for  $\alpha\in\delta$, 
and $Y_\emptyset=\emptyset$.
Then
$(\pkl-\widehat{2})\cup\build{\nabla}_{\alpha<\delta}^{} 
X_\alpha=(\pkl-\widehat{2})\cup
\build{\nabla}_{e\in P_2(\delta)}^{}Y_e$.

\bigskip
\item[2)] Let $X_e\subseteq\pkl$ for $e\in P_2(\delta)$. Define 
$Y_\alpha$ for $\alpha<\delta$ by
$Y_\alpha=X_{\{\alpha\}}$. Then

\centerline{$(\pkl-\widehat{2})\cup X_\emptyset ~\cup
\build{\nabla}_{\alpha<\delta}^{} Y_\alpha=(\pkl-\widehat{2}) ~\cup
\build{\nabla}_{e\in P_2(\delta)}^{} X_e$. }
\end{enumerate}
\vskip-0,5cm\eop

We finally turn to the question of existence of $\dt$-normal ideals. 
Let us first deal with the
degenerate case $\kappa=\omega$.\\

\hskip-0,2cm\pmb{\sc Proposition 2.10.} \  {\it Assume 
$\kappa=\omega$. Then there exists a
$\dt$-normal ideal on $\pkl$ if  and only if $\delta<\omega$.}\\
\vskip-0,5cm
\vspace*{0ex}{\it Proof. }The right-to-left implication is
immediate from Proposition 2.2 $(1)$. For the \linebreak reverse
implication, observe that
$\pol=(\pol-\widehat{2})~\cup\build{\nabla}_{e\in P_2(\omega)}^{} B_e$, where
$B_\emptyset=\emptyset$  and $B_{\{n\}}=\{a\in\pol:\cup(a\cap\omega)=n\}$ for
$n\in\omega$. Hence $\pol\in\nabla^{
\lbrack\omega\rbrack^{<2}}\iol$ by Lemma 2.1 $((0)$ and $(1))$. If 
$\delta\ge\omega$, then
$\pol\in\nabla^{\dt}J$ by Lemma 2.1 $(2)$, and therefore $J$ is not 
$\dt$-normal.\eop
\vskip-0,3cm
We will now look for sufficient conditions for the existence of 
$\dt$-normal ideals on $\pkl$
in the case $\kappa>\omega$. We will use the following key lemma.\\

\hskip-0,2cm\pmb{\sc Lemma 2.11.} \ {\it
\vskip-0,3cm
\begin{enumerate}\vskip-0,3cm
\item[$(0)$] Assume $\theta\cdot\aleph_0<\kappa$ and 
$|P_{\theta}(\mu)|<\kappa$ for every cardinal
$\mu<
\kappa$. Then $\pkl\not\in\nabla^{\lbrack\lambda\rbrack^{<\theta}}\ikl$.
\item[$(1)$] {\rm ([M1])} Assume that $\kappa$ is Mahlo. Then
$\pkl\not\in\nabla^\lk\ikl$.
\end{enumerate}}

\vspace*{0ex}{\it Proof.}
\vskip-0,3cm
\begin{enumerate}\vskip-0,3cm
\item[(0)] : Let $b_e\in\pkl$ for $e\in
P_\theta(\lambda)$, and  fix $a\in\pkl$. Set $\rho = 
\theta\cdot\aleph_0$ if $\theta\cdot\aleph_0$ is
regular, and
$\rho=(\theta\cdot\aleph_0)^+$ otherwise.  Now define $x_\alpha\in\pkl$ for
$\alpha<\rho$ so that
\item[(i)]\vskip-0,2cm $x_0=a\cup\theta$.
\item[(ii)] \vskip-0,2cm If $\alpha>0$, then 
$\displaystyle\bigcup_{\beta<\alpha}x_\beta\subseteq
x_\alpha$ and $x_\alpha\in
\displaystyle\bigcap \ \{\widehat{b_e}:e\in 
P_\theta(\displaystyle\bigcup_{\beta<\alpha}x_\beta)\}$.
\end{enumerate}

Set $x=\displaystyle\bigcup_{\alpha<\rho}x_\alpha$.
  Given $e\in P_{|x\cap\theta|}(x)$, there is $\beta<\rho$ with $e\in 
P_\theta(x_\beta)$. Then
$b_e\subseteq x_{\beta+1}\subseteq x$. Thus
$\widehat{a}~\cap\build{\Delta}_{e\in{P}_\theta(\lambda)}^{}\widehat{b_e}
\not=\emptyset$. Hence 
$\pkl\not\in\nabla^{\lbrack\lambda\rbrack^{<\theta}}\ikl$
  by Proposition 2.3 ((iv) $\rightarrow$ (i)).
\bigskip
\begin{enumerate}
\item[(1)] : Let $b_e\in\pkl$ for $e\in\pkl$, and fix $a\in\pkl$. 
Define $x_\alpha\in\pkl$ and
$\gamma_\alpha\in\kappa$ for $\alpha<\kappa$ \ so that

\item[(i)]\vskip-0,2cm $\gamma_\alpha=\cup(x_\alpha\cap \kappa)$.
\item[(ii)]\vskip-0,2cm$x_0=a$.
\item[(iii)]\vskip-0,2cm $x_{\alpha+1}\in\widehat{x_\alpha}\cap
\widehat{\mbox{ }\{\gamma_\alpha+1\} \mbox{ }}\cap\displaystyle\bigcap_{e
\subseteq x_\alpha}\widehat{b_e}$.
\item[(iv)]\vskip-0,2cm 
$x_\alpha=\displaystyle\bigcup_{\beta<\alpha}x_\beta$ if $\alpha$ is 
an infinite
limit ordinal  $>0$.
\end{enumerate}

There is a regular infinite cardinal $\tau$ such that $\gamma_\tau=\tau$. Then
$x_\tau\in\widehat{a}\ \cap\build{\Delta}_{e\in\pkl}^{}\widehat{b_e}$. Hence
$\pkl\not\in\nabla^\lk\ikl$  by Proposition 2.3 ((iv) $\rightarrow$ (i)).\eop

{\bf Definition.} \ {\it For  $f:\ptd\longrightarrow\pkl$, 
$C_f^{\kappa,\lambda}$ denotes the
set  of all $a\in\pkl$  such that $a\cap\theta\not=\emptyset$ and 
$f(e)\subseteq a$ for every $e\in
P_{|a\cap\theta|}(a\cap \delta)$.}

\bigskip
The following is straightforward.\\

\hskip-0,2cm\pmb{\sc Lemma 2.12.} \  {\it Given $B\subseteq\pkl$, 
$B\in\nabla^\dt\ikl$ if and
only if $B\cap C_f^{\kappa,\lambda}=
\emptyset$ for some $f:\ptd\longrightarrow\pkl$.}\\

\hskip-0,2cm\pmb{\sc Lemma 2.13.} \  {\it Assume that 
$\delta\ge\kappa$ and either
$\theta=\kappa$ and
$\kappa$ is Mahlo, or
  $3\le\theta$, $\theta\cdot\aleph_0<\kappa$ and 
$|P_\theta(\mu)|<\kappa$ for every cardinal
$\mu<\kappa$. Then $\nabla^\dt\ikl$ is a
$\dt$-normal ideal on $\pkl$.}
\bigskip

\vspace*{0ex}{\it Proof. }$\nabla^\dt\ikl$ is an ideal on $\pkl$ by
Lemma 2.11, Lemma 2.1 $(2)$ and Proposition 2.3
((i) $\rightarrow$ (ii)).

\bigskip
Assume  $\theta\ge\omega$. Given $g_b:\ptd\longrightarrow\pkl$ for 
$b\in\ptd$, define $f:\ptd
\longrightarrow\pkl$ by
$f(e)=\displaystyle\bigcup_{b,c\ \in P_\theta(e)}g_b(c)$. Then 
$\widehat{\omega}\cap
C_f^{\kappa,\lambda}\subseteq
\build{\Delta}_{b\in\ptd}^{} C_{g_b}^{\kappa,\lambda}$. Hence 
$\nabla^\dt\ikl$ is
\mbox{$\dt$-normal} by  Lemma 2.12.

\bigskip
Now assume  $3\le\theta<\omega$. Select a bijection 
$j:\ptd\longrightarrow P_2(\delta)$.
Given
$g_b:\ptd\longrightarrow\pkl$  for $b\in\ptd$, define 
$f:\ptd\longrightarrow\pkl$ by

\centerline{$f(e)=\displaystyle\bigcup~\{g_b(c):b,c\in\ptd$ and 
$j(b)\cup j(c)\subseteq
  e\}$.}

Then $\widehat{\theta}\cap C_j^{\kappa,\lambda}\cap
C_f^{\kappa,\lambda}\subseteq\build{\Delta}_{b\in\ptd}^{} 
C_{g_b}^{\kappa,\lambda}$.
Hence $\nabla^\dt\ikl$ is $\dt$-normal by Lemma 2.12.\eop

\bigskip
\hskip-0,2cm\pmb{\sc Lemma 2.14.} \ {\it Assume  $J$ is
$\lbrack\delta\rbrack^{<\overline{\theta}\cdot 3}$-normal. Then $J$ is
$\dt$-normal.}\bigskip

\vspace*{2ex}{\it Proof. }If $\overline{\theta}\ge 3$,
$J=\nabla^{\lbrack\delta\rbrack^{<
\overline{\theta}\cdot 3}}
J=\nabla^\dt J$ by Proposition 2.2 (0). If $\overline{\theta}<3$, \break
$J\subseteq \nabla^\dt J
\subseteq\nabla^{\lbrack\delta\rbrack^{<3}}J\subseteq J$ by
Lemma 2.1 ((0) and (2)).\eop

It remains to show that our sufficient conditions are also necessary ones.\\

\bigskip
\hskip-0,2cm\pmb{\sc Lemma 2.15.} \ {\it Assume 
$\pkl\not\in\nabla^\dt\ikl$, and let
$\mu,\tau$  be two cardinals such that $\mu<\kappa\cap(\delta+1)$ and
$0<\tau<\theta^+\cap\kappa$.  Then $|P_\tau(\mu)|<\kappa$.}

\bigskip
\vspace*{0ex}{\it Proof. }Suppose otherwise, and pick a one-to-one $j:
\kappa\longrightarrow P_\tau(\mu)$. Define\break  $f:
\widehat{\mu\cup\tau}\longrightarrow P_\tau(\mu)$ by 
$f(a)=j(\cup(a\cap\kappa))$.
Then  $f$ is $\ptd$-regressive, which contradicts \mbox{Proposition 2.3}
((i) $\rightarrow$ (v)).\eop

\hskip-0,2cm\pmb{\sc Lemma 2.16.}
\begin{enumerate}
\item[(0)]\vskip-0,2cm {\it  Assume that $\delta\geq\kappa>\omega$ 
and $\delta$ is a limit ordinal.
Then}
\vskip-1cm
$$\{a\in\pkl:\cup(a\cap\delta) \ \hbox{\rm is a limit ordinal and} \
\cup(a\cap\delta)\not\in 
a\cap\delta\}\in(\nabla^{\lbrack\delta\rbrack^{<2}}\ikl)^*.$$

\item[(1)] {\it Assume $\delta\geq\kappa > \omega$. Then the set of 
all $a\in P_\kappa(\lambda)$
such that $cf (\cup (a\cap\eta)) < |a\cap\overline\theta|$ for some 
limit ordinal $\eta$ with
$\kappa \leq \eta \leq \delta$ and $cf (\eta ) \geq \overline\theta$ lies in
$\nabla^{[\delta]^{{<\theta}}} I_{\kappa,\lambda}$.}
\item[(2)] {\it Assume  $\kappa>\omega$, and let $C$ be a closed 
unbounded subset of $\kappa$.
Then}

\centerline{$\{a\in
\pkl:a\cap\kappa\in C\}\in(\nabla^{\lbrack\kappa\rbrack^{<2}}
\ikl)^*$. }
\end{enumerate}

\vspace*{0ex}{\it Proof. }Use Lemma 2.11 (0) and Proposition 2.3.\eop

\hskip-0,2cm\pmb{\sc Lemma 2.17.} \  {\it Assume that $\kappa$ is an 
uncountable  limit
cardinal  and
$\pkl\not\in\nabla^{\lbrack\kappa
\rbrack^{<\kappa}}\ikl$. Then $\kappa$ is Mahlo.}

\vspace*{2ex}{\it Proof. }By Lemma 2.1 (3) and Lemma 2.16.\eop

Our study of the case $\kappa>\omega$ culminates in the following\\

\hskip-0,2cm\pmb{\sc Proposition 2.18.} \
\begin{enumerate}
\item[(0)] {\it  Assume that $\kappa>\omega$. Further assume that 
$\delta<\kappa$, or
$\theta<\kappa$, or $\kappa$ is  not a limit cardinal. Then there 
exists a $\dt$-normal ideal on
$\pkl$ if and only if
$|P_{\overline{\theta}}(\mu)| <\kappa$ for every cardinal 
$\mu<\kappa\cap(\delta+1)$.}

\item[(1)] {\it Assume that $\delta\ge\kappa>\omega$, $\theta=\kappa$ 
and $\kappa$ is a limit
cardinal. Then there exists a \break
  $\dt$-normal ideal on $\pkl$ if and only if $\kappa$ is Mahlo.}
\end{enumerate}

\vspace*{0ex}{\it Proof. }
\begin{enumerate}
\item[(0)] : Let us first assume that there exists a $\dt$-normal ideal on
$\pkl$. Then $\pkl\not\in
\nabla^\dt\ikl$ by Lemma 2.1 (2).
Notice that if $\delta<\kappa$, $\theta=\kappa$ and $\kappa$ is a 
limit cardinal, then
setting
$\tau=|\delta|^+$,  we have that $\tau<\theta^+\cap\kappa$ and
$P_{\overline{\theta}}(|\delta|)=P_\tau(|\delta|)$.  Hence by Lemma 2.15,
$|P_{\overline{\theta}}(\mu)|<\kappa$ for every cardinal 
$\mu<\kappa\cap(\delta+1)$.
\bigskip

Conversely, assume that $|P_{\overline{\theta}}(\mu)|<\kappa$ for 
every cardinal $\mu<\kappa\cap(\delta+1)$.
If $\delta<\kappa$, then $|P_{\overline{\theta}\cdot 
3}(\delta)|<\kappa$, and therefore
$\ikl$ is
$\lbrack\delta\rbrack^{<\overline{\theta}\cdot 3}$-normal by 
Proposition 2.2 (1).
If $\delta\ge\kappa$, then $\overline{\theta}<\kappa$, and 
consequently $\nabla^{
\lbrack\delta\rbrack^{<\overline{\theta}\cdot 3}}\ikl$ is a
$\lbrack\delta\rbrack^{<\overline{\theta}\cdot 3}$-normal
ideal on $\pkl$ by Lemma 2.13. Thus by Lemma 2.14 there exists a 
$\dt$-normal ideal on $\pkl$.
\item[(1)] : If $\kappa$ is Mahlo, then $\nabla^\dt\ikl$ is a 
$\dt$-normal ideal on $\pkl$ by Lemma 2.13.
Conversely, if there exists a $\dt$-normal ideal on $\pkl$, then 
$\kappa$ is Mahlo by Lemma 2.1
(2) and Lemma 2.17.\eop
\end{enumerate}

\hskip-0,2cm\pmb{\sc Corollary 2.19.} \  {\it There exists a 
$\dt$-normal ideal on $\pkl$ if
and only if there exists a
$\lbrack
\delta\cap\kappa\rbrack^{<\theta\cap|\delta|^+}$-normal ideal on 
$P_\kappa(\kappa)$.}

\vspace*{2ex}{\it Proof. }By Propositions 2.10 and 2.18.\eop

\hskip-0,2cm\pmb{\sc Corollary 2.20.} \  {\it Assume that  $\delta < 
\kappa$ and there exists
a
$[\delta]^{<\theta}$ - normal ideal on $P_\kappa(\lambda)$. Then every ideal on
$P_\kappa(\lambda)$ is $[\delta]^{<\theta}$ - normal.}

\vspace*{2ex}{\it Proof. }By Propositions 2.2 and 2.18(0).\eop

The following (see e.g. [EHM\' aR]) is due independently to Hajnal 
and Shelah.\\

\hskip-0,2cm\pmb{\sc Lemma 2.21.} \  {\it  Let $\mu$ be an infinite 
cardinal. Then
$\mu^\rho$  assumes only finitely many values for $\rho$ with $2^\rho < \mu$.
}\\

\hskip-0,2cm\pmb{\sc Lemma 2.22.} \  {\it  Let $\mu, \chi$  be two infinite cardinals suchthat \  $2^{<\chi} \leq\mu$. Then  ${(\mu^{<\chi})}^ 
{<\chi} = \mu^{<\chi}$. }
\\
\vspace*{0ex}{\it Proof. } If there exists a cardinal $\tau < \chi$ 
such that $2^\tau = \mu$,
then
\vskip-1cm
$$\mu^{<\chi} = {(2^\tau)}^{<\chi} = 2^{<\chi} = \mu.$$
Otherwise, there exists by Lemma 2.21 a cardinal $\rho < \chi$ such 
that $\mu^{<\chi} =
\mu^\rho$. Then
$${(\mu^{<\chi})}^{<\chi} = {(\mu^\rho)}^{<\chi} = \mu^{<\chi}. $$
\vskip-1cm\eop

\hskip-0,2cm\pmb{\sc Proposition 2.23.} \  {\it Assume that there exists a
$[\kappa]^{<\theta}$ - normal ideal on $P_\kappa(\lambda)$. Then} (a) ([M3])
$\kappa^{<\overline\theta} = \kappa$, {\it and} (b)
${(\mu^{<\overline\theta})}^{<\overline\theta} = 
\mu^{<\overline\theta}$ {\it  for every
cardinal $\mu > \kappa$.}
\\

\vspace*{0ex}{\it Proof. }  (b) follows from Lemma 2.22 since by 
Proposition 2.18 \
$2^{<\overline\theta} \leq \kappa$.\eop

\vskip1cm

{\Large\bf 3. $NS_{\kappa,\lambda}^\dt$}\\

In this section we  describe the smallest $\dt$-normal ideal on 
$\pkl$. We will need the
following which shows that we could without loss of generality assume 
$\theta$ to be an
infinite cardinal.\\

\hskip-0,2cm\pmb{\sc Lemma 3.1.} \  {\it Assume $J$ is $\dt$-normal. 
Then $J$ is
$\lbrack\delta\rbrack^{<\theta\cdot \aleph_0}$-normal.}

\vspace*{2ex}{\it Proof. } We can assume that $\theta<\omega$ since 
otherwise the result is
trivial. The desired conclusion is immediate from  Proposition 2.2 
(1) in case $\delta<\omega$.
Now assume  $\delta\ge\omega$. We have  $\kappa>\omega$ by
Proposition  2.10. Fix $A\in J^+$ and a $P_\omega(\delta)$-regressive 
$f:A\longrightarrow
P_\omega(\delta)$. We  define a $\ptd$-regressive 
$g:A\cap\widehat{\omega}\longrightarrow\ptd$
by
$g(a)=\{|f(a)|\}$.  By Proposition 2.5 \break ((i) $\rightarrow$ 
(iv)), there are $C\in J^+\cap
P(A\cap\widehat{\omega})$ and $n\in\omega$  such that $g$ is 
identically $n$ on $C$. If
$n=0$,  $f$ is clearly constant on $C$. Otherwise,  select a bijection
$j_a:n\longrightarrow f(a)$ for each $a\in C$. Using Proposition 2.5 
((i)$\rightarrow$
(iv)), define
$C_k\in J^+$ for $k\le n$ and $h_i:C_i\longrightarrow\ptd$ for $i<n$ so that
\begin{enumerate}
\item[(0)]\vskip-0,2cm $C_0=C$.
\item[(1)]\vskip-0,2cm $C_{i+1}\subseteq C_i$.
\item[(2)]\vskip-0,2cm $h_i(a)=\{j_a(i)\}$.
\item[(3)]\vskip-0,2cm $h_i$ is constant on $C_{i+1}$.
\end{enumerate}
Then $f$ is constant on $C_n$. Hence $J$ is 
$\lbrack\delta\rbrack^{<\omega}$-normal by
Proposition 2.5
((iv) $\rightarrow$ (i)).\eop

\bigskip
\hskip-0,2cm\pmb{\sc Proposition 3.2.} \  {\it If there exists a 
$\dt$-normal ideal on
$\pkl$,  then the smallest such ideal is
$\nabla^{\lbrack\delta\rbrack^{<\overline{\theta}\cdot 3}}\ikl$.}
\bigskip

\vspace*{0ex}{\it Proof. }Assume that there exists a $\dt$-normal 
ideal on $\pkl$. Then
$\nabla^{\lbrack\delta\rbrack^{<\overline{\theta}\cdot 
3}}\ikl\subseteq K$ for every
$\dt$-normal ideal $K$  on $\pkl$
  by Lemmas 3.1 and  2.1 (2). Moreover,
$\nabla^{\lbrack\delta\rbrack^{<\overline{\theta}\cdot 3}}\ikl$  is 
itself a $\dt$-normal ideal on
$\pkl$ by the proofs of Propositions 2.10 and 2.18.\eop

{\bf Definition.}  {\it Assuming the existence of a $\dt$-normal 
ideal on $\pkl$, we set
\break
$NS_{\kappa,\lambda}^\dt=\nabla^{\lbrack\delta\rbrack^{<\overline{\theta}\cdot
3}}\ikl$.}\\

\hskip-0,2cm\pmb{\sc Proposition 3.3.} \  {\it Let $\delta^\prime$  be an
ordinal with  $1\leq \delta'\leq\delta$,
and $\theta'$ be a cardinal with $2\leq\theta'\le\theta$. Then $
NS_{\kappa,\lambda}^{\lbrack\delta^\prime
\rbrack^{<\theta^\prime}} \subseteq\nsdt$.}

\vspace*{2ex}{\it Proof. }By Proposition 2.4.\eop

\hskip-0,2cm\pmb{\sc Proposition 3.4.} \  {\it
$\nsdt=NS_{\kappa,\lambda}^{\lbrack\delta\rbrack^{<\theta\cdot\aleph_0}}
=NS_{\kappa,\lambda}^{\lbrack\delta\rbrack^{<\overline{\theta}}}$.}
\bigskip

\vspace*{0ex}{\it Proof. } We have
$NS_{\kappa,\lambda}^{\lbrack\delta\rbrack^{<\theta\cdot\aleph_0}}\subseteq\nsdt\subseteq N 
S_{\kappa,\lambda}^{\lbrack\delta\rbrack^{<\overline{\theta}}}\subseteq 
NS_{\kappa,\lambda}^{\lbrack
\delta\rbrack^{<{\theta}\cdot\aleph_0}}$ by Lemma 3.1 and 
Propositions 3.2, 2.2 (0) and 3.3.\eop

\bigskip
\hskip-0,2cm\pmb{\sc Proposition 3.5.} \  {\it $\nsdt=\ikl$ \ if \ 
$\delta<\kappa$.} \\
\vspace*{2ex}{\it Proof. }By Corollary 2.20.\eop

{\bf Definition.} {\it We put
$NS_{\kappa,\lambda}^\delta=NS_{\kappa,\lambda}^{\lbrack\delta\rbrack^{<2}}$ .}

\bigskip
It follows from Propositions 2.9 and 3.2 that 
$NS_{\kappa,\lambda}^\delta$ is the smallest
$\delta$-normal ideal on
$\pkl$. We will conform to usage and denote $NS_{\kappa,\lambda}^\lambda$ by
$NS_{\kappa,\lambda}$. \\

The following is due to Abe [A].

\bigskip
\hskip-0,2cm\pmb{\sc Proposition 3.6.} \  {\it Assume 
$\kappa\le\delta<\kappa^+$. Then
$NS_{\kappa,\lambda}^\delta=\nabla^{\lbrack\delta
\rbrack^{<2}}\ikl$.}\\
\vskip-0,8cm
\vspace*{2ex}{\it Proof. }Let us first prove the assertion for
$\delta=\kappa$. Given $f_b:P_2(\kappa)
\longrightarrow\pkl$ for
$b\in P_2(\kappa)$, define $f:P_2(\kappa)\longrightarrow\pkl$ by
$f(e)=\displaystyle\bigcup_{b\in P_2((\cup e)+1)}
\ \displaystyle\bigcup_{c\in P_2((\cup e)+1)}f_b(c)$.
Then  $C_f^{\kappa,\lambda}
\subseteq$

$\build{\Delta}_{b\in P_2(\kappa)}^{} C_{f_b}^{\kappa,\lambda}$.
Hence $\nabla^{\lbrack\kappa\rbrack^{<2}}\ikl$ is 
$\lbrack\kappa\rbrack^{<2}$-normal by
Lemma 2.12 and
Proposition 2.5 ((ii) $\rightarrow$ (i)). It follows that
$NS_{\kappa,\lambda}^\kappa=\nabla^{\lbrack\kappa\rbrack^{ <2}}\ikl$ 
by Proposition 3.2 and
Lemma 2.1 (2).

\bigskip
Now assume  $\kappa<\delta<\kappa^+$. By Propositions 2.18 (0) and 
4.4 (below), there
exists \break $A\in(
\nabla^{\lbrack
\delta\rbrack^{<2}}\ikl)^*$ such that 
$NS_{\kappa,\lambda}^{\lbrack\delta\rbrack^{<2}}=
NS_{\kappa,\lambda}^{
\lbrack\kappa\rbrack^{<2}}|A$. Then by Lemma 2.1 (2),

\centerline{\hspace*{3.3cm}$\nabla^{\lbrack\delta\rbrack^{<2}}\ikl\subseteq
NS_{\kappa,\lambda}^{\lbrack\delta\rbrack^{<2}}=(\nabla^{\lbrack\kappa 
\rbrack^{<2}}\ikl)|A\subseteq
  \nabla^{\lbrack\delta\rbrack^{<2}}\ikl$.\hspace{3.4cm}$\boite$}

\bigskip\bigskip
Abe [A] also showed that for \ $\delta\geq \kappa^+$, 
$NS_{\kappa,\lambda}^\delta -
\nabla^{[\delta]^{<2}} I_{\kappa,\lambda} \not= \phi$ (in fact

$\nabla^{[\kappa]^{<2}}
(\nabla^{[\kappa^+]^{<2}} I_{\kappa,\lambda}) - 
\nabla^{[\kappa]^{<2}} I_{\kappa,\lambda} \not=
\phi)$.

\bigskip
By Lemma 2.12, $NS_{\kappa,\lambda}^{[\delta]^{<\theta}}$ is the set of all
$B\subseteq P_\kappa (\lambda)$ such that $B\cap C_f^{\kappa,\lambda} 
= \phi$ for some $f :
P_{\overline\theta\cdot 3} (\delta) \rightarrow P_\kappa(\lambda)$. 
The following generalizes a
well-known (see Lemma 1.13 in [Me] and Proposition 1.4 in [M2]) 
characterization of
$NS_{\kappa,\lambda}$.\\

\bigskip
\hskip-0,2cm\pmb{\sc Proposition 3.7.} \  {\it Assume 
$\delta\ge\kappa$. Then given
$B\subseteq\pkl$, $B\in\nsdt$ if  and only if $B\cap\{a\in
C_g^{\kappa,\lambda}:a\cap\kappa\in\kappa\}=\emptyset$ for some
$g:P_{\overline{\theta}.3}(\delta)\longrightarrow P_3(\lambda)$.}\\

\vskip-0,2cm
\vspace*{0ex}{\it Proof. }Set $\tau=2$ if $\overline{\theta}<\omega$
and $\delta<\kappa^+$, $\tau=3$ if
$\overline{\theta}<\omega$ and $\delta\ge\kappa^+$, and
$\tau=\overline{\theta}$
if
$\overline{\theta}\ge\omega$.
Then by Lemma 2.12 and Propositions 3.4 and 3.6, it suffices
to show that for every
$f:P_\tau(\delta)\longrightarrow\pkl$, there exists
$g:P_{\overline{\theta}\cdot 3}(\delta)
\longrightarrow P_3(\lambda)$ with the
property that $\{a\in C_g^{\kappa,\lambda}:a\cap\kappa\in\kappa\}
\subseteq C_f^{\kappa,\lambda}$.
Thus fix $f:P_\tau(\delta)\longrightarrow
\pkl$. Pick a bijection $j_e:|f(e)|\longrightarrow f(e)$ for each
$e\in P_\tau(\delta)$.

\bigskip
Let us first assume that $\overline{\theta}\ge\omega$. Define 
$h:P_\tau(\delta)\longrightarrow\kappa$
by

\centerline{$h(e)=\omega\cup(((\cup(e\cap\kappa))+1)+|f(e)|)$.}

We define
$k:P_\tau(\delta)\longrightarrow\lambda$ as follows.
Given $e\in P_\tau(\delta)$, set $\alpha=\cup(e\cap\kappa)$. We put $k(e)=0$
if
$\alpha\not\in e$. Assuming now that $\alpha\in e$, put 
$c=e-\{\alpha\}$ and $\xi=\cup(c\cap\kappa)$,
and let $\beta$ denote the unique ordinal $\zeta$ such that 
$\alpha=(\xi+1)+\zeta$. We put
$k(e)=j_c(\beta)$ if $\beta\in|f(c)|$ and $k(e)=0$ otherwise. Finally  define
$g:P_\tau(\delta)\longrightarrow P_3(\lambda)$ by  $g(e)=\{h(e),k(e)\}$.
Now fix $a\in C_g^{\kappa,\lambda}$ with $a\cap\kappa\in\kappa$, and 
$c\in P_{|a\cap\tau|}(a\cap\delta)$. Put $\xi=\cup(c\cap\kappa)$.
Given $\beta\in|f(c)|$, set $e=c\cup\{(\xi+1)+\beta\}$. Since 
$h(c)\subseteq a$, we have
$\omega\subseteq a$ and $(\xi+1)+\beta\in a$, and therefore $e\in 
P_{|a\cap\tau|}(a\cap\delta)$.
Hence $j_c(\beta)\in a$, since clearly $k(e)=j_c(\beta)$. Thus 
$f(c)\subseteq a$.

\bigskip
Let us next assume that $\overline{\theta}<\omega$ and 
$\delta\ge\kappa^+$. Select a bijection
$h:P_3(\delta)\longrightarrow\delta-\kappa$. Define $k:P_3(\delta)
\longrightarrow\lambda$ so that (a) $k(\emptyset)=2$, and (b) given 
$e\in P_3(\delta)$,
$k(\{h(e)\})=|f(e)|$
and for all $\beta\in|f(e)|$, $k(\{\beta,h(e)\})=j_e(\beta)$. Then define
  $g:P_{\overline{\theta}\cdot 3}(\delta)\longrightarrow P_3(\lambda)$ so that
$g(e)=\{h(e),k(e)\}$ for all $e\in P_3(\delta)$. It is readily 
checked  that $g$
is as desired.

\eject
Finally, assume that $\overline{\theta}<\omega$ and $\delta<\kappa^+$.
Define $h:P_2(\delta)\longrightarrow\kappa$ by~:

\begin{enumerate}
\item[(i)] $h(\emptyset)=2+|f(\emptyset)|$.

\item[(ii)] $h(\{\alpha\})=(\alpha+1)+|f(\{\alpha\})|$ for  $\alpha\in\kappa$.

\item[(iii)] $h(\{\alpha\})=|f(\{\alpha\})|$ for  $\alpha\in\delta-\kappa$.
\end{enumerate}
Then define $k:P_3(\delta)\longrightarrow\lambda$ so that
\begin{enumerate}
\item[(0)] $k(\{\beta\})=j_\emptyset(\beta)$ whenever $\beta\in|f(\emptyset)|$.

\item[(1)] $k(\{\alpha,(\alpha+1)+\beta\})=j_{\{\alpha\}}(\beta)$ whenever
$\alpha\in\kappa$ and $\beta\in|f(\{\alpha\})|$.

\item[(2)] $k(\{\alpha,\beta\})=j_{\{\alpha\}}(\beta)$ whenever
$\alpha\in\delta-\kappa$ and
$\beta\in|f(\{\alpha\})|$.\\
  \end{enumerate}

Finally define $g:P_{\overline{\theta}\cdot 3}(\delta)\longrightarrow
P_3(\lambda)$ so that
$g(e)=\{h(e),k(e)\}$ if $e\in P_2(\delta)$, and $g(e)=\{k(e)\}$ if $e\in
P_3(\delta)-P_2(\delta)$. Then
  $g$ is as desired.\eop

\vskip1cm

{\Large\bf 4. Variations of $\delta$}\\

This section is concerned with the case when $\delta$ is not a cardinal.

{\bf Throughout the section it is assumed that $\delta\geq\kappa$}.

Our first remark is that   we do not lose generality by assuming that
$\delta =\kappa\alpha$ for some ordinal $\alpha>0$. Lemma 4.1 and 
Proposition 4.2
generalize results of
\mbox{Abe [A].}\\

\hskip-0,2cm\pmb{\sc Lemma 4.1.} \  {\it Assume that 
$\delta=\kappa\alpha$ for some ordinal
$\alpha>0$, and $J$ is
$\dt$-normal.  Then $J$ is $\lbrack\delta+\xi\rbrack^{<\theta}$
-normal for every $\xi<\kappa$.}\\
\vspace*{2ex}{\it Proof. }Fix $\xi<\kappa$. Since
$\xi+\kappa\alpha=\kappa\alpha$, we can define
$j:\kappa\alpha+\xi
\longrightarrow\kappa\alpha$ by :
  $j(\beta)=\xi+\beta$ for  $\beta<\kappa\alpha$, and 
$j(\kappa\alpha+\gamma)=\gamma$
for  $\gamma<\xi$. Set

\centerline{$C=\widehat{\xi}\cap\{a\in\pkl:(\forall\beta\in 
a\cap\kappa\alpha)\;j(\beta)\in
a\}$.}
\bigskip
Then clearly $C\in(\nsdt)^*$. Now given $A\in J^+$ and a 
$P_\theta(\delta+\xi)$-regressive
\break $f:A
\longrightarrow
  P_\theta(\delta+\xi)$, define $g:A\cap C\longrightarrow 
P_\theta(\delta)$ by  $g(a)=j\lbrack
  f(a)\rbrack$. Since $A\cap C\in J^+$ by Proposition 3.2, and $g$ is 
$\ptd$-regressive, we have
by  Proposition 2.5 ((i)$\rightarrow$ (iv)) that $g$ is constant on 
some $D\in J^+$. Then $f$
is constant  on $D$.
Hence $J$ is $\lbrack\delta+\xi\rbrack^{<\theta}$-normal by 
Proposition 2.5 ((iv)
$\rightarrow$ (i)).\eop

\vskip-0,5cm
\hskip-0,2cm\pmb{\sc Proposition 4.2.} \  {\it Assume that 
$\delta=\kappa\alpha$ for some
ordinal $\alpha>0$. Then
\vskip-0,3cm
\begin{enumerate}
\vskip-0,3cm
\item[{\rm (a)}] 
$\nsdt=NS_{\kappa,\lambda}^{\lbrack\delta+\xi\rbrack^{<\theta}}$ for 
every \
$\xi<\kappa$.
\item[{\rm (b)}] $NS_{\kappa,\lambda}^{\lbrack\delta+\kappa\rbrack^{<2}}-\nsdt
\not=\emptyset$.
\end{enumerate}}

\vspace*{0ex}{\it Proof. }\vskip-0,5cm
\begin{enumerate}
\vskip-0,7cm
\item[(a)] : By Lemma 4.1, Propositions 3.2 and 3.3.
\vskip-0,7cm
\item[(b)] : Select 
$f:\lbrack\delta+\kappa\rbrack^{<3}\longrightarrow\pkl$ so that
$f(\{\beta\})=\{\beta +1\}$ for every $\beta\in\delta+\kappa$. Given $g:
P_{\overline\theta\cdot 3}(\delta)
\longrightarrow\pkl$, pick $a\in C_g^{\kappa,\lambda}$ and 
$\gamma\in(\delta+\kappa)-\delta$ with
$\gamma\ge\cup(a\cap(\delta+\kappa))$. Then $a\cup\{\gamma\}\in 
C_g^{\kappa,\lambda}-C_f^{\kappa,
\lambda}$.
Hence $\pkl-C_f^{\kappa,\lambda}\in 
NS_{\kappa,\lambda}^{\lbrack\delta+\kappa\rbrack^{<2}}-\nsdt$ by
Lemma 2.12. \end{enumerate}\eop
 
\vskip-1,3cm
\hskip-0,2cm\pmb{\sc Lemma 4.3.} \  {\it  The following are equivalent~:
\vskip-0,3cm
\begin{enumerate}
\vskip-0,3cm
\item[{\rm (i)}] $J$ is $[\delta]^{<\theta}$- normal.
\item[{\rm (i)}] $\nabla^\delta I_{\kappa,\lambda} \subseteq J$ and 
$J$ is $[|\delta|]^{<\theta}$-
normal.
\end{enumerate}}

\vspace*{0ex}{\it Proof. }  (i) $\rightarrow$ (ii)~: By Lemma 2.1 (2).

\hskip1,4cm (ii) $\rightarrow$ (i)~: Select a bijection $j : \delta 
\rightarrow |\delta|$
and set
$$D = P_\kappa(\lambda) - \build{\nabla}_{\alpha < \delta}^{} \big( 
P_\kappa(\lambda) -
\widehat{\{ j(\alpha)\}}\big).$$
Then $D$ lies in $(\nabla^\delta I_{\kappa,\lambda})^*$ and so in 
$J^*$. Now fix $A\in J^+$ and
a $P_{\overline\theta\cdot 3} (\delta)$-regressive \break $f : A\longrightarrow
P_{\overline\theta\cdot 3} (\delta)$. Define $g : A\cap D \longrightarrow
P_{\overline\theta\cdot 3} \  (|(\delta|)$ by $g(a) = j[f(a)]$. Since $g$ is
$P_{\overline\theta.3} (|\delta|)$-regressive,  we can find $C\in J^+ 
\cap P(A\cap D)$ and $u
\in P_{\overline\theta\cdot 3} (|\delta|)$ so that $g(a) = u$ for all 
$a\in C$. Then $f$ takes
the constant value $j^{-1} (u)$ on $C$.
\eop

Let us remark in passing that Lemma 4.3 can be combined with a result 
of [M4] to show that $J$
is $[\delta]^{<\theta}$-normal if and only if it is $\delta$-normal and $(\mu,
|\delta|)$-distributive for every infinite cardinal $\mu < \overline\theta$.

\bigskip
\hskip-0,2cm\pmb{\sc Proposition 4.4.} \  {\it 
$NS_{\kappa,\lambda}^{[\delta]^{<\theta}} =
NS_{\kappa,\lambda}^{[|\delta|]^{<\theta}}\mid D$  for some $D\in 
(\nabla^\delta
I_{\kappa,\lambda})^\ast$.}
\bigskip

\vspace*{0ex}{\it Proof. } By the proof of Lemma 4.3. \eop

Using Cantor's normal form for the base $|\delta|$, one easily 
obtains the following.

\bigskip
\hskip-0,2cm\pmb{\sc Proposition 4.5.} \  {\it Assume that $\gamma < 
\delta \leq
\gamma^\gamma$, where
$\gamma = |\delta|$. Then $NS_{\kappa,\lambda}^\delta = 
NS_{\kappa,\lambda}^\gamma \mid A$,
where $A$ is the set of all $a\in P_\kappa (\lambda)$ with the 
following property~: Suppose
that $1\leq\alpha < \delta$ and $\alpha = \gamma^{\eta_1} \xi_1 
+\ldots + \gamma^{\eta_p}
\xi_{p}$, where $1\leq p < \omega$, $\gamma > \eta_1 >\ldots > 
\eta_p$ and $\gamma > \xi_i \geq
1$ for $1\leq i\leq p$. Then $\alpha\in a$ if and only if $\{ \eta_1, 
\xi_1,\ldots,\eta_p,
\xi_p\} \subseteq a$.}

\bigskip
Thus for example $NS_{\kappa,\lambda}^{\kappa+\kappa} = 
NS_{\kappa,\lambda}^\kappa \mid A$,
where $A$ is the set of all $a\in P_\kappa (\lambda)$ such that
$$a-\kappa = \{ \kappa+\alpha : \alpha\in a\cap \kappa\},$$
and $NS_{\kappa,\lambda}^{\kappa^2} = NS_{\kappa,\lambda}^\kappa \mid 
B$, where $B$ is the
set of all $a\in P_\kappa (\lambda)$ such that
$$a-\kappa = \{ \kappa\beta+\alpha : \alpha,\beta\in a\cap\kappa \ 
\hbox{and} \ \beta\geq
1\}.$$

\vskip1,5cm


{\Large\bf 5. Variations of $\theta$}\\

\hskip-0,2cm\pmb{\sc Proposition 5.1.} \  {\it Assume that 
$\delta\ge\kappa$ and
$\overline{\theta}\cdot\aleph_0$ is a regular cardinal, and let
$\theta^\prime$ be a cardinal  such that $\theta'\leq\kappa$ and
$\overline{\theta}\cdot\aleph_0<\overline{\theta^\prime}$.  Then
$NS_{\kappa,\lambda}^{[\kappa]^{<\theta'}}
-\nsdt\not=\emptyset$  (and therefore 
$NS_{\kappa,\lambda}^{[\delta]^{<\theta}} \not=
NS_{\kappa,\lambda}^{[\delta]^{<\theta'}})$.}

\vspace*{2ex}{\it Proof. } Given $f:P_{\overline{\theta}\cdot
3}(\delta)\longrightarrow\pkl$, we use  Proposition 2.18 (0) to 
define $a_\alpha\in\pkl$ and
$\gamma_\alpha\in\kappa$ for $\alpha <\overline{\theta}\cdot\aleph_0$ 
as follows:
\begin{enumerate}
\item[(i)]\vskip-0,2cm $a_0=\overline{\theta}\cdot 3$.

\item[(ii)] \vskip-0,2cm $\gamma_\alpha=\cup(a_\alpha\cap\kappa)$.

\item[(iii)] \vskip-0,2cm 
$a_{\alpha+1}=a_\alpha\cup(\gamma_\alpha+1)\cup(\cup f\lbrack
P_{\overline{\theta}\cdot 3} (a_\alpha\cap\delta)\rbrack)$.

\item[(iv)] \vskip-0,2cm 
$a_\alpha=\displaystyle\bigcup_{\beta<\alpha}a_\beta$ if $\alpha$ is 
an infinite
limit ordinal.
\end{enumerate}

Put 
$a=\displaystyle\bigcup_{\alpha<\overline{\theta}\cdot\aleph_0}a_\alpha$. 
Then
$a\in C_f^{\kappa,\lambda}$ and
$cf(\cup(a\cap\kappa))=\overline{\theta}\cdot\aleph_0$. Hence

\centerline{$\{a\in\pkl: cf(\cup(a\cap\kappa))=\overline{\theta}\cdot\aleph_0
\}\in(\nsdt)^+$}

by Lemma 2.12. It remains to observe that

\centerline{$\{a\in\pkl:cf(\cup(a\cap\kappa))>\overline{\theta}\cdot\aleph_0\}\in
(\nabla^{\lbrack\kappa\rbrack^{<\overline{\theta^\prime}}}\ikl)^*$}

by Lemma 2.16 (1).\eop

\bigskip
We will see that the conclusion of Proposition 5.1 may fail if
$\overline\theta\cdot\aleph_0$ (and hence $\theta$) is a singular 
cardinal. The remainder of
the section is concerned with the case when $\theta$ is a limit cardinal.

The following is immediate from Proposition 2.18 (0).

\bigskip
\hskip-0,2cm\pmb{\sc Proposition 5.2.} \  {\it  Suppose that $\theta$ 
is a limit cardinal
with $\theta < \kappa$. Then the following are equivalent~:}

\begin{enumerate}
\item[(i)]\vskip-0,3cm {\it   There exists a 
$[\delta]^{<\theta}$-normal ideal on
$P_\kappa(\lambda)$.}
\item[(ii)] \vskip-0,5cm {\it For each cardinal $\rho$ with 
$2\leq\rho<\theta$, there exists a
$[\delta]^{<\rho}$-normal ideal on $P_\kappa(\lambda)$.}
\end{enumerate}

Notice that if $\theta = \kappa$ and $\kappa$ is an inacessible 
cardinal that is not Mahlo,
then by Proposition 2.18, (ii) holds but (i) does not.\\

\hskip-0,2cm\pmb{\sc Proposition 5.3.} \  {\it Assume that 
$\delta\ge\kappa$ and $\theta$ is
a limit cardinal. Then the following are equivalent~:}
\begin{enumerate}
\item[(i)] \vskip-0,5cm {\it $J$ is $[\delta]^{<\theta}$-normal.}

\item[(ii)] \vskip-0,5cm {\it $J$ is $[\delta]^{<\rho}$-normal for 
every cardinal $\rho$ with
$2\leq\rho < \theta$.}
\end{enumerate}

\vspace*{2ex}{\it Proof.} \ (i) $\rightarrow$ (ii) : By Lemma 2.1 (2).

\hskip1,2cm (ii) $\rightarrow$ (i) : By Proposition 2.5, it suffices 
to show that if $A\in
J^+$ and $f : A\rightarrow P_\theta (\delta)$ is $P_\theta(\delta)$-regressive,
then $f\upharpoonright D$ is $P_\rho (\delta)$-regressive for some 
$D\in J^+\cap
P(A)$ and some cardinal $\rho$ with $2\le \rho < \theta$. This is clear if
$\theta < \kappa$. Assuming $\theta = \kappa$, put $B = \{ a\in A : a
\cap\kappa\in\kappa\}$. Then $\mid f(a)\mid \in a\cap\kappa$ for every $a\in B$
with $a\cap\kappa\not= \emptyset$. It remains to observe that by Lemmas 2.1 (2)
and 2.16 (2), $J$ is $[\kappa]^{<2}$-normal and $B\cap\widehat 2\in J^+$.
\eop

We have the following corresponding characterization of
$NS_{\kappa,\lambda}^{[\delta]^{<\theta}}$. \\

\hskip-0,2cm\pmb{\sc Proposition 5.4.} \  {\it  Assume that 
$\delta\ge\kappa$ and $\theta$ is
a limit cardinal. Then \quad
\break
$NS_{\kappa,\lambda}^{[\delta]^{<\theta}} = \nabla^\theta
(\displaystyle\bigcup_{2\le\rho <\theta}
NS_{\kappa,\lambda}^{[\delta]^{<\rho}})$.
}\\

\vspace*{2ex}{\it Proof.} We have $\nabla^\theta
(\displaystyle\bigcup_{2\le\rho<\theta} NS_{\kappa,\lambda}^{[\delta]^{<\rho}})
\subseteq \nabla^\theta NS_{\kappa,\lambda}^{[\delta]^{<\theta}} \subseteq
NS_{\kappa,\lambda}^{[\delta]^{<\theta}}$ by Proposition 2.9 and Lemma 2.1 (2).
Now fix $B\in NS_{\kappa,\lambda}^{[\delta]^{<\theta}}$. Then by Lemma 2.12,
there is $f : P_\theta(\delta) \longrightarrow P_\kappa(\lambda)$ such that
$B\cap C_f^{\kappa,\lambda} = \emptyset$. Set $f_\rho = f\upharpoonright P_\rho
(\delta)$ for each cardinal $\rho$ with $2\le\rho < \theta$. Let us define
$D\subseteq P_\kappa(\lambda)$ by~:  $D = \widehat\theta$ if $\theta < \kappa$,
and
$$D = \{ a\in P_\kappa(\lambda) : a\cap\kappa \ \hbox{ is an infinite 
limit cardinal}
\}$$
otherwise. Then $D\in {(NS_{\kappa,\lambda}^{[\delta]^{<2}})}^*$ by Lemmas 2.16
(2) and 2.1 (2). Let $A$ be the set of all $a\in D$ such that $a\in
C_{f_{\rho}}^{\kappa,\lambda}$ for every cardinal $\rho$ with 
$2\le\rho < \theta
\cap (a\cap\kappa)$. Then clearly, $A\subseteq C_f^{\kappa,\lambda}$ and
$P_\kappa(\lambda) - A\in\nabla^\theta 
(\displaystyle\bigcup_{2\le\rho < \theta}
NS_{\kappa,\lambda}^{[\delta]^{<\rho}}).$ Hence, $B\in\nabla^\theta
(\displaystyle\bigcup_{2\le\rho <\theta}
NS_{\kappa,\lambda}^{[\delta]^{<\rho}})$.
\eop

Now we focus on the case when $\theta$ is a singular cardinal. \\

\hskip-0,2cm\pmb{\sc Proposition 5.5.} \  {\it Assume that there exists a
$[\delta]^{<\theta}$-normal ideal on $P_\kappa(\lambda)$, $\theta$ is 
a singular cardinal,
and either $\delta\geq 2^{<\theta}$, or $\delta\geq\theta$ and $cf 
(\theta^{<\theta}) \not=
cf (\theta)$. Then there exists a $[\delta]^{<\theta^+}$-normal ideal 
on $P_\kappa
(\lambda)$.}
\bigskip

\vspace*{0ex}{\it Proof.}  Notice that by Proposition 2.18 (0), 
$2^{<\theta} \leq
\theta^{<\theta} < \kappa$. First suppose that $\theta\leq\delta < 
2^{<\theta}$ and $cf
(\theta^{<\theta}) \not= cf(\theta)$. Then there is a cardinal $\tau 
< \theta$ such that
$\theta^{<\theta} = \theta^\tau$. We get
$$|\delta|^\theta  \leq {(2^\theta)}^\theta = \theta^\theta = 
{(\theta^{<\theta})}^{cf
(\theta)} = \theta^{\tau \cdot cf(\theta)} = \theta^{<\theta},$$
so the desired conclusion follows from Proposition 2.18 (0). Now 
suppose $\delta\geq
2^{<\theta}$. Let $\mu$ be a cardinal with $2^{<\theta} \leq \mu < 
\kappa \cap (\delta
+1)$. Then by Lemma 2.22 and Proposition 2.18 (0),
$$\mu^\theta = {(\mu^{<\theta})}^{<\theta} = \mu^{<\theta} < \kappa.$$
 From this together with Proposition 2.18 (0), we get the desired conclusion.
\eop

Observe that if $\theta$ is a singular cardinal with $cf 
(\theta^{<\theta}) = cf (\theta)$,
then for $\delta = \theta$ and \break $\kappa = 
{(\theta^{<\theta})}^+$, (a) there is a
$[\delta]^{<\theta}$-normal ideal on $P_\kappa (\lambda)$, but (b) there is no
$[\delta]^{<\theta^+}$-normal ideal on $P_\kappa (\lambda)$ \ (since 
$\theta^\theta =
{(\theta^{<\theta})}^{cf (\theta)} \geq \kappa$). \\

Let us recall the following definition.

\bigskip
{\bf Definition.} \ {\it Given cardinals $\sigma,\nu,\rho$ with 
$2\le\sigma$ and $\sigma\cdot
\aleph_0 \leq\nu$, ${\rm cov} (\rho, \nu, \nu, \sigma)$ is the least 
cardinal $\mu$
for which one can find $X\subseteq P_\nu(\rho)$ such that $\mid X\mid 
= \mu$ and
for every $c\in P_\nu(\rho)$, there is $d\in P_\sigma(X)$ with $c\subseteq \cup
d$. }

\bigskip
\hskip-0,2cm\pmb{\sc Lemma 5.6.} ([S1])\  {\it Let $\sigma, \nu, 
\rho$  be three cardinals
such that \
$2\leq\sigma$ and $\sigma\cdot\aleph_0 \leq\nu < \rho$.  Then the 
following hold~:
\begin{enumerate}
\item[\rm (i)]\vskip-0,2cm ${\rm cov}~(\rho, \nu, \nu, \sigma) \geq \rho$.

\item[\rm (ii)]\vskip-0,3cm ${\rm cov}~(\rho^+, \nu, \nu, \sigma) = 
\rho^+ \cdot {\rm cov} (\rho, \nu,
\nu, \sigma)$.

\item[\rm (iii)]\vskip-0,3cm If $cf (\rho) < \sigma$, then ${\rm 
cov}~(\rho, \nu,\nu,\sigma) =
\displaystyle\bigcup_{\nu < \mu < \rho} {\rm cov}~(\mu, \nu,\nu,\sigma)$.
\end{enumerate}
}

\bigskip
\hskip-0,2cm\pmb{\sc Proposition 5.7.} \  {\it  Assume that $\theta$ is a
singular cardinal, $\delta\geq\kappa$ and there is a cardinal
$\sigma$ such that $2\leq\sigma < \theta$ and  $\rm{cov}(\mid\delta\mid, \theta,\theta, \sigma) = \ \mid\delta\mid$. Then every $[\delta]^{<\ 
theta}$-normal
ideal on $P_\kappa (\lambda)$ is $[\delta]^{<\theta^+}$-normal.} \\
\vspace*{0ex}{\it Proof.} Assume $J$ is $[\delta]^{<\theta}$-normal. 
Since by Proposition
2.18 (0) $2^{<\theta} < \delta$, we can find
$x_\xi\in P_\theta (\delta)$ for $\xi\in\delta$, and $f : 
P_\theta(\delta) \longrightarrow
P_\sigma(\delta)$ so that $c =
\displaystyle\bigcup_{\xi\in f(c)} x_\xi$  for every $c\in P_\theta(\delta)$.
Now fix  $A_e\in J^*$ for $e\in P_{\theta^+} (\delta)$. Put $B_d =
A_{\build{\cup}_{\xi\in D}^{}x_\xi}$ for $d\in P_\theta(\delta)$. Set $C =
\build{\Delta}_{c\in P_{\theta}(\delta)}^{} \widehat{f(c)}$, $D =
\build{\Delta}_{d\in P_{\theta}(\delta)}^{} B_d$ and $E = C \cap D \cap
\widehat\theta$. Then $E\in J^*$ by Proposition 2.5. Let $a\in E$ and $e\in
P_{\mid a\cap\theta^+\mid} (a\cap\delta)$ be given. Select $c_\zeta\in
P_\theta(\delta)$ for $\zeta < cf(\theta)$ so that $e =
\displaystyle\bigcup_{\zeta < cf (\theta)} c_\zeta$. For each $\zeta < cf
(\theta)$, we get $c_\zeta\in P_{\mid a\cap\theta\mid} (a\cap\delta)$ and
therefore $f(c_\zeta) \subseteq a$. So setting $d =
\displaystyle\bigcup_{\zeta<cf (\theta)} f(c_\zeta)$, we have $d\in P_{\mid
a\cap\theta\mid} (a\cap\delta)$ and consequently $a\in B_d$. Notice that $B_d =
A_e$, since
$$\bigcup_{\xi\in d} x_\xi = \bigcup_{\zeta < cf(\theta)} \bigcup_{\xi\in
f(c_\zeta)} x_\xi = \bigcup_{\zeta < cf(\theta)} c_\zeta = e.$$
Thus $E\subseteq \build{\Delta}_{e\in P_{\theta^{+}}(\delta)}^{} A_e$, and
therefore $\build{\Delta}_{e\in P_{\theta^+}(\delta)}^{} A_e \in J^*$. Hence by
Proposition 2.5, $J$ is $[\delta]^{<\theta^+}$-normal. \eop

\hskip-0,2cm\pmb{\sc Corollary 5.8.} \  {\it  Assume that $\theta$ is 
a singular cardinal and
$\delta\geq\kappa$. Assume further that either $\theta$ is a strong 
limit cardinal, or
$\delta < \kappa^{+\theta}$. Then every $[\delta]^{<\theta}$-normal ideal on
$P_\kappa(\lambda)$ is $[\delta]^{<\theta^+}$-normal (and hence
$NS_{\kappa,\lambda}^{[\delta]^{<\theta^+}}= 
NS_{\kappa,\lambda}^{[\delta]^{<\theta}})$.
}
\bigskip

\vspace*{0ex}{\it Proof. } The result follows from Proposition 5.7 
and the following
remarks~: (a) If $\theta$ is a singular strong limit cardinal, then 
by a result of Shelah
[S2], for every cardinal $\rho > \theta$, there is a cardinal 
$\sigma$ such that
$2\leq\sigma < \theta$ and cov $(\rho, \theta,\theta,\sigma) = \rho$. 
(b) If $n < \omega$,
then by Proposition 2.23
\quad cov $(\kappa^{+n}, \theta,\theta,2) = {(\kappa^{+n})}^{<\theta}
= \kappa^{+n}.$ \quad
(c) Using Lemma 5.6, it is easy to show by induction that if
$\omega\leq\gamma < \theta$,  then cov$(\kappa^{+\gamma}, 
\theta,\theta, |\gamma|^+) =
\kappa^{+\gamma}$. \eop

\vskip1,5cm

{\Large\bf 6. The case $\kappa\leq\delta < \kappa^{+\overline\theta}$}\\

{\bf Definition.} \ {\it $E_{\kappa,\lambda}$ denotes the set of all 
$a\in P_\kappa(\lambda)$
such that $a\cap\kappa \not= \phi$ and $a\cap\kappa = \cup (a\cap\kappa)$.}

\eject
\hskip-0,2cm\pmb{\sc Proposition 6.1.} ([M3]) \ {\it Assuming the 
existence of a
$[\kappa]^{<\theta}$-normal ideal on $P_\kappa (\lambda)$, the 
following are equivalent~:

\begin{enumerate}
\item[\rm (i)]\vskip-0,3cm $J$ is $[\kappa]^{<\theta}$-normal.

\item[\rm (ii)]\vskip-0,3cm $J$ is $\kappa$-normal and  \
$\{ a\in E_{\kappa,\lambda} : cf (a\cap\kappa) \geq \cup (a\cap 
\overline\theta)\} \in
J^*$.
\end{enumerate}}

We will show that this result can be generalized.

\bigskip
{\bf Definition.} \ {\it Let $\rho$ be a cardinal with 
$\kappa\leq\rho$, and $\beta$ be an
ordinal with $1\leq\beta < \kappa$. Then 
$A_{\kappa,\lambda}^{\rho,\beta}$ denotes the set
of all $a\in P_\kappa(\lambda)$ such that (a) $\alpha + 1\in a$ for 
every $\alpha\in a\cap
(\rho^{+\beta} - \rho)$, and (b) $\rho^{+\gamma}\in a$ for every 
$\gamma < \beta$.}

\bigskip
Thus if $a\in A_{\kappa,\lambda}^{\rho,\beta}$ and $\gamma < \beta$, 
then $\cup (a\cap
\rho^{+(\gamma +1)})$ is a limit ordinal that is strictly greater 
than $\rho^{+\gamma}$ and
does not belong to a.

\bigskip
\hskip-0,2cm\pmb{\sc Proposition 6.2.} \ {\it Assume that $\delta = 
\rho^{+\beta}$, where
$\rho$ is a cardinal with $\kappa\leq\rho$, and $\beta$ an ordinal 
with $1\leq\beta <
\overline\theta$.  Then the following are equivalent~:}
\begin{enumerate}
\item[(i)] {\it $J$ is $[\delta]^{<\theta}$-normal.}

\item[(ii)]\vskip-0,2cm {\it $J$ is $[\delta]^{<|\beta|^+}$-normal 
and $[\rho]^{<\theta}$-normal, and
the set of all
$a\in A_{\kappa,\lambda}^{\rho,\beta}$ such that \break $cf (\cup 
(a\cap \rho^{+(\alpha+1)}))
\geq \cup (a\cap\overline\theta)$ for every $\alpha <\beta$ lies in $J^*$.}
\end{enumerate}

\vspace*{0ex}{\it Proof. } (i) $\rightarrow$ (ii) : By Lemma 2.16 (1) 
and Propositions 2.4
and 6.1.

\hskip1,3cm (ii) $\rightarrow$ (i)~: By Proposition 5.3 it suffices 
to prove the result for
$\theta < \kappa$. We can also assume that $|\beta|^+ < \theta$ 
(since otherwise the
result is trivial) and (by Lemma 3.1) that $\theta$ is an infinite cardinal.

For $\gamma\in\delta - \rho$, select a bijection $\widetilde\gamma : 
\gamma \longrightarrow
|\gamma|$. Let $B$ be the set of all $a\in 
A_{\kappa,\lambda}^{\rho,\beta}$ such that
($\ast$)
$\theta \subseteq a$, ($\ast\ast$) $cf (\cup (a \cap \rho^{+(\alpha 
+1)}))\geq \theta$  for
all
$\alpha < \beta$, and ($\ast\ast\ast$) $\widetilde\gamma (\xi) \in a$ 
whenever $\gamma\in
a \cap (\delta - \rho)$ and $\xi\in a\cap\gamma$.  Notice that $B\in 
J^*$. For $a\in B$ and
$\alpha < \beta$, select $z_\alpha^a \subseteq a \cap (\kappa^{+(\alpha+1)} -
\kappa^{+\alpha})$ so that  o.t. $(z_\alpha^a) = cf (\cup (a\cap 
\kappa^{+(\alpha+1)}))$
and $\cup z_\alpha^a = \cup (a\cap \kappa^{+(\alpha+1)})$.

Now fix $C\in J^+$ and a $P_\theta (\delta)$-regressive $F : 
C\longrightarrow P_\theta
(\delta)$. Set $D = C\cap B$. For $a\in D$ and $1\leq\eta\leq\beta$, 
define $k_\eta^a :
P_\theta (a\cap\rho^{+\eta}) \longrightarrow P_{|\eta|^+} 
(a\cap\rho^{+\eta})$ as follows~:

\begin{enumerate}
\item[{\rm (0)}] $k_1^a (e) = \{ \gamma\}$, where $\gamma$ is the 
least $\zeta\in
z_0^a$ such that $e\subseteq\zeta$.

\item[{\rm (1)}] If $e-\rho^{+\eta} \not= \phi$, $k_{\eta+1}^a (e) = 
\{ \gamma\} \cup
k_\eta^a (\widetilde\gamma [e])$,  where $\gamma$ is the least 
$\zeta\in z_\eta^a$ such that
$e\subseteq\zeta$. Otherwise $k_{\eta+1}^a (e) = k_\xi^a (e)$,  where 
$\xi$ is the least
$\chi \geq 1$ such that $e\subseteq \rho^{+\chi}$.

\item[{\rm (2)}] Suppose that $\eta$ is a limit ordinal. If $\cup e = 
\kappa^{+\eta}$,
$k_\eta^a (e) = \displaystyle\bigcup_{\alpha < \eta} k_{\alpha+1}^a (e \cap
\rho^{+(\alpha+1)})$. Otherwise $k_\eta^a (e) = k_\xi^a (e)$, where 
$\xi$ is the least
$\chi\geq 1$ such that $e\subseteq \rho^{+\chi}$.
\end{enumerate}

\bigskip
Let $a\in D$. For $1\leq\xi\leq\beta$, let $\Phi_\xi$   assert that 
given $\zeta\in e\in
P_\theta (a\cap \rho^{+\xi})$, there are $n\in\omega$ and 
$\gamma_0,\ldots,\gamma_n\in
k_\xi^a (e)$ such that $\zeta \in\gamma_0$, $(\widetilde\gamma_j 
\circ \ldots \circ
\widetilde\gamma_0) (\zeta) \in \gamma_{j+1}$  for $j=0,\ldots,n-1$, and
$(\widetilde\gamma_n \circ \ldots \circ \widetilde\gamma_0) (\zeta) 
\in a\cap\rho$.  Let
us show by induction that $\Phi_\xi$ holds. Given $\zeta\in e\in 
P_\theta (a\cap\rho^+)$,
let $k_1^a (e) = \{ \gamma\}$.  Then $e\subseteq\gamma$ and 
$\widetilde\gamma (\zeta)\in
a\cap\rho$. Thus $\Phi_1$ holds. Next suppose that 
$1<\alpha\leq\beta$ and $\Phi_\xi$
holds for $1\leq\xi <\alpha$. Let $\zeta\in e\in P_\theta (a\cap 
\rho^{+\alpha})$, where
$e-\rho^{+\xi} \not= \phi$ for every $\xi < \alpha$. Define $\xi, 
\gamma_0$ and $e'$ as
follows~:
\begin{enumerate}
\item[(a)] If  $\alpha$ is a limit ordinal, $\xi$ is the least 
$\sigma$ such that
$\zeta\in\rho^{+(\sigma+1)}$.  Otherwise $\xi+1 = \alpha$.
\item[(b)] $\gamma_0\in z_\xi^a \cap k_{\xi+1}^a (e\cap \rho^{+(\xi +1)})$.
\item[(c)] $e\cap \rho^{+(\xi+1)} \subseteq \gamma_0$.
\item[(d)] $e' = \widetilde\gamma_0 [e\cap\rho^{+(\xi+1)}]$
\item[(e)] $k_\xi^a (e') \subseteq k_{\xi+1}^a (e\cap \rho^{+(\xi+1)})$.
\end{enumerate}

Then $\xi < \alpha$  and $\zeta\in\gamma_0\in z_\xi^a \cap k_\alpha^a 
(e)$. Moreover
$\widetilde\gamma_0(\zeta) \in e'\in P_\theta (a\cap \rho^{+\xi})$ 
and \break $k_\xi (e')
\subseteq k_\alpha^a (e)$. If $\xi = 0$, then 
$\widetilde\gamma_0(\zeta)\in a\cap\rho$.
Otherwise, there are $\gamma_1,\ldots,\gamma_n\in k_\xi^a (e')$, 
where $1\leq n <
\omega$, such that  $\widetilde\gamma_0(\zeta)\in\gamma_1$, 
$(\widetilde\gamma_j
\circ\ldots\circ\widetilde\gamma_1) 
(\widetilde\gamma_0(\zeta))\in\gamma_{j+1}$ for
$j=1,\ldots,n-1$, and $(\widetilde\gamma_n\circ\ldots\circ\widetilde\gamma_1)
(\widetilde\gamma_0 (\zeta))\in a\cap\rho$. So $\Phi_\alpha$ holds.

\bigskip
Define $G : D \longrightarrow P_{|\beta|^+} (\delta)$ by $G(a) = 
k_\beta^a \big(F(a)\big)$.
Since $G$ is $P_{|\beta|^+} (\delta)$-regressive, there are $T\in J^+\cap
P(D)$ and $x\in P_{|\beta|^+} (\delta)$ such that $G$ takes the 
constant value $x$ on $T$.
For $a\in T$ and $\zeta\in F(a)$, pick $\chi_\zeta^a \in a\cap\rho$ 
so that there exist
$n\in\omega$ and $\gamma_0,\ldots,\gamma_n\in x$ such that $\zeta\in\gamma_0$,
$(\widetilde\gamma_j \circ\ldots\circ\widetilde\gamma_0) 
(\zeta)\in\gamma_{j+1}$ for
$j=0,\ldots,n-1$, and 
$(\widetilde\gamma_n\circ\ldots\circ\widetilde\gamma_0)(\zeta) =
\chi_\zeta^a$.  Now define $H : T\longrightarrow P_\theta (\rho)$ by $H(a) = \{
\chi_\zeta^a : \zeta\in F(a)\}$. Since $H$ is $P_\theta 
(\rho)$-regressive, we can find
$W\in J^+ \cap P(T)$ and $y\in P_\theta (\rho)$ so that $H$ takes the 
constant value $y$ on
$W$. Let $d$ be the set of all $\zeta\in\delta$ for which one can 
find $n\in\omega$ and
$\gamma_0,\ldots,\gamma_n\in x$ so that $\zeta\in\gamma_0$, 
$(\widetilde\gamma_j
\circ\ldots\circ\widetilde\gamma_0) (\zeta) \in \gamma_{j+1}$ for 
$j=0,\ldots,n-1$ and
$(\widetilde\gamma_n\circ\ldots\circ\widetilde\gamma_0)(\zeta)\in y$. 
Then $|d| < \theta$
and $F[W] \subseteq P_\theta (d)$. Since $|P_\theta (d)| < \kappa$ by 
Proposition 2.18
(0), there are $Z\in J^+ \cap P(W)$ and $v\in P_\theta (d)$ such that 
$F$ takes the
constant value $v$ on $Z$.  \eop

\hskip-0,2cm\pmb{\sc Corollary 6.3.} \ {\it Assume that $|\delta| = 
\kappa^{+n}$, where
$n<\omega$.  Then $NS_{\kappa,\lambda}^{[\delta]^{<\theta}} = 
NS_{\kappa,\lambda}^\delta |
C$,  where
$C$ is the set of all $a\in P_\kappa (\lambda)$ such that $cf (\cup 
(a\cap \kappa^{+m})) \geq
\cup (a\cap \overline\theta)$ for every $m\leq n$.
}\\
\vspace*{0ex}{\it Proof. } By Lemma 4.3 and Propositions 6.1 and 6.2.  \eop

\hskip-0,2cm\pmb{\sc Corollary 6.4.} \ {\it  Assume that $|\delta| = 
\kappa^{+\beta}$, where
$\omega
\leq\beta < \theta$. Then $NS_{\kappa,\lambda}^{[\delta]^{<\theta}} =
NS_{\kappa,\lambda}^{[\delta]^{<|\beta|^+}}\mid~C$, where $C$ is the 
set of all $a\in
P_\kappa(\lambda)$ such that {\rm (a)} $cf (\cup (a\cap\kappa))\geq 
\cup (a\cap\overline\theta)$,
and {\rm (b)} $cf\big(\cup (a\cap\kappa^{+(\alpha+1)})\big)\geq \cup 
(a\cap\overline\theta)$ for
every $\alpha < \beta$.
}
\bigskip

\vspace*{0ex}{\it Proof. } By Lemma 4.3 and Propositions 6.1 and 6.2. \eop

So for example for $\kappa > \omega_2$ and $\lambda = \kappa^{+\omega}$,
$NS_{\kappa,\lambda}^{[\lambda]^{<\aleph_2}} =
NS_{\kappa,\lambda}^{[\lambda]^{<\aleph_1}}\mid C$,  where $C$ is the 
set of all $a\in
P_\kappa(\lambda)$ such that $cf (\cup (a\cap\kappa^{+n})) \geq 
\omega_2$ for every
$n<\omega$. We will see later (see Corollary 9.6) that if 
$\lambda^{\aleph_0} = 2^\lambda$, then
$NS_{\kappa,\lambda}^{[\lambda]^{<\aleph_2}}\mid A \not=
NS_{\kappa,\lambda}^{[\lambda]^{<\aleph_0}}\mid A$ for all $A$.

\vskip1,5cm


{\Large\bf 7. $NS_{\kappa,\lambda}^{[\delta]^{<\theta}} | A$}
\bigskip

In this section we continue to investigate whether given 
$\delta'\geq\delta$ and
$\theta'\geq\theta$ with \break $(\delta',\theta') \not= 
(\delta,\theta)$, it is possible to find
$A$ such that $NS_{\kappa,\lambda}^{[\delta']^{<\theta'}} =
NS_{\kappa,\lambda}^{[\delta]^{<\theta}}\mid A$. The \break following 
is obvious.\\

\eject
\hskip-0,2cm\pmb{\sc Lemma 7.1.} \ {\it Let $\delta'$ be an ordinal with
$\delta\leq\delta'\leq\lambda$, and $\theta'$ be a cardinal with
$\theta\leq\theta'\leq\kappa$.  Then the  following are equivalent~:

\begin{enumerate}
\item[\rm (i)] There exists $A\in {\big( 
NS_{\kappa,\lambda}^{[\delta]^{<\theta}}\big)}^+$
such that \  $NS_{\kappa,\lambda}^{[\delta']^{<\theta'}} =
NS_{\kappa,\lambda}^{[\delta]^{<\theta}}\mid A$.
\item[\rm (ii)] There is  $f : P_{\overline{\theta'}\cdot 3} 
(\delta') \longrightarrow
P_\kappa(\lambda)$ such that  for every  $h : 
P_{\overline{\theta'}\cdot 3} (\delta')
\longrightarrow P_\kappa(\lambda)$, one can find $k : 
P_{\overline\theta\cdot 3} (\delta)
\longrightarrow P_\kappa(\lambda)$ with $C_f^{\kappa,\lambda} \cap 
C_k^{\kappa,\lambda}
\subseteq C_h^{\kappa,\lambda}$.
\end{enumerate}}

We start with a positive result.

\bigskip
\hskip-0,2cm\pmb{\sc  Lemma 7.2.} \  {\it Let $\delta'$ be an ordinal with
$\delta\leq\delta'\leq\lambda$, and $\theta'$ be a cardinal with
$\theta\leq\theta'\leq\kappa$.\break  Assume that $\delta\geq\kappa$ and
$|\delta|^{<\overline\theta} = |\delta'|^{<\overline{\theta'}}$. Then
$NS_{\kappa,\lambda}^{\lbrack\delta^\prime\rbrack^{<\theta^\prime}}=\nsdt
|A$ for some
\break  $A\in(\nabla^{
\lbrack\delta^\prime\rbrack^{<\overline{\theta^\prime}}}\ikl)^*$.}
\bigskip

\vspace*{0ex}{\it Proof. }Select a bijection $j: 
P_{\overline{\theta'}} (\delta')
\longrightarrow P_{\overline\theta} (\delta)$ with $j(\phi) = \phi$,
and let $i$ denote its inverse.  Define
$f: P_{\overline{\theta'}} (\delta')\longrightarrow\pkl$  by~:
$f(b)=(\overline{\theta}\cdot 3)\cup j(b)$ if
$\overline{\theta}<\kappa$, and
$f(b)=|j(b)|^+\cup j(b)$ otherwise. Then
$C_f^{\kappa,\lambda}\in(\nabla^{\lbrack\delta^\prime\rbrack^{
<\overline{\theta^\prime}}}\ikl)^*$ by Lemma 2.12. Now given
$h:P_{\overline{\theta^\prime}\cdot 3}(\delta')
\longrightarrow\pkl$, define
$k:P_{\overline{\theta}\cdot 3}(\delta)\longrightarrow\pkl$ so that

\begin{enumerate}
\item[(i)] $k(e)=(h\circ i)(e)$ whenever $e\in P_{\overline{\theta}}(\delta)$;

\item[(ii)] If $\overline{\theta^\prime}=2$, then 
$k(\{\alpha,\beta\})=h(i(\{\alpha\})\cup
i(\{\beta\}))$ whenever
$\alpha$ and $\beta$ are two distinct members of $\delta$.
\end{enumerate}
It is readily checked that $C_f^{\kappa,\lambda}\cap 
C_k^{\kappa,\lambda}\subseteq
C_h^{\kappa,\lambda}$.  Hence
$NS_{\kappa,\lambda}^{\lbrack\delta^\prime\rbrack^{<\theta^\prime}}=
\nsdt|C_f^{\kappa,
\lambda}$.\eop

\hskip-0,2cm\pmb{\sc  Lemma 7.3.} \  {\it Assume that there exists a
$[\kappa]^{<\theta}$-normal ideal on $P_\kappa(\lambda)$. Let $\nu$ 
be a cardinal with
$\nu > \kappa$, and $\sigma$ be the least cardinal $\tau$ with 
$\tau^{<\overline\theta}\geq
\nu$.  Then {\rm (a)} $\sigma >
\kappa$, {\rm (b)} $\mu^{<\overline\theta} < \sigma$ for every 
cardinal $\mu < \sigma$, {\rm (c)}
$\sigma^{<\overline\theta} = \nu^{<\overline\theta}$, and {\rm (d)} 
$\sigma^{<\overline\theta} =
\sigma$ if $cf(\sigma) \geq \overline\theta$, and $\sigma^{<\overline\theta} =
\sigma^{cf(\sigma)}$ otherwise.
}\\

\vspace*{0ex}{\it Proof.} Proposition 2.23 tells us that \ 
$\kappa^{<\overline\theta} =
\kappa$,\  so $\sigma > \kappa$. Moreover given a  \break cardinal 
$\mu$  with $\kappa < \mu <
\sigma$, we have $\mu^{<\overline\theta} < \sigma$ since otherwise by 
Proposition 2.23
\break $\mu^{<\overline\theta} = 
{(\mu^{<\overline\theta})}^{<\overline\theta} \geq
\sigma^{<\overline\theta} \geq \nu$,  which would contradict the 
definition of $\sigma$.
Again by Proposition 2.23, $\sigma^{<\overline\theta} =
{(\sigma^{<\overline\theta})}^{<\overline\theta} \geq \nu^{<\overline\theta}$
and hence $\sigma^{<\overline\theta} = \nu^{<\overline\theta}$. 
Finally, for (d) use
e.g. Lemma 1.7.3 in [HoStW]. \eop

\hskip-0,2cm\pmb{\sc  Proposition 7.4.} \  {\it Assume 
$\delta\geq\kappa$, and let $\sigma$
the least cardinal $\tau$ such that $\tau^{<\overline\theta} \geq 
|\delta|$. \break Then\
$NS_{\kappa,\lambda}^{[\delta]^{<\theta}} = 
NS_{\kappa,\lambda}^\sigma |A$ for some
$A\in {\big( \nabla^{[\delta]^{<\overline\theta}} 
I_{\kappa,\lambda}\big)}^*$ if
$cf(\sigma) \geq \overline\theta$, and \break
$NS_{\kappa,\lambda}^{[\delta]^{<\theta}} =
NS_{\kappa,\lambda}^{[\sigma]^{<(cf(\sigma))^+}} | D$ for some $D \in
{\big(\nabla^{[\delta]^{<\overline\theta}} 
I_{\kappa,\lambda}\big)}^*$  otherwise. }
\bigskip

\vspace*{0ex}{\it Proof. } By Lemmas 7.2 and 7.3.~\eop
Lemma 7.2 has the following generalization.

\bigskip
\hskip-0,2cm\pmb{\sc  Proposition 7.5.} \  {\it Assume 
$|\delta'|^{<\overline{\theta'}} =
|\delta|^{<\overline\theta}$, where  $\delta'$ is an ordinal with
$\kappa\leq\delta'\leq\lambda$, and $\theta'$ a cardinal with
$2\leq\theta'\leq\kappa$. Then $NS_{\kappa,\lambda}^{[\delta']^{<\theta'}}| C =
NS_{\kappa,\lambda}^{[\delta]^{<\theta}} | C$ for some
  $C\in {\big(
\nabla^{[\delta^{\prime\prime}]^{<\overline{\theta^{\prime\prime}}}}
I_{\kappa,\lambda}\big) }^*$,  where $\delta^{\prime\prime} = 
\delta\cup\delta'$
and $\theta^{\prime\prime} = \theta\cup\theta'$.
}
\bigskip

\vspace*{0ex}{\it Proof. } By Lemma 7.2 we can find $A,B\in {\big(
\nabla^{[\delta^{\prime\prime}]^{<\overline{\theta^{\prime\prime}}}}I_ 
{\kappa,\lambda}\big)}^*$
so that $NS_{\kappa,\lambda}^{[\delta]^{<\theta}} | A =
NS_{\kappa,\lambda}^{[\delta^{\prime\prime}]^{<\theta^{\prime\prime}}} =
NS_{\kappa,\lambda}^{[\delta']^{<\theta'}} | B$. Then $C = A \cap B$ 
is as desired.\eop

We will now describe some situations when $\delta\leq\delta'$, $\theta
\leq\theta'$, $|\delta|^{\overline\theta} < 
|\delta'|^{<\overline{\theta'}}$  and there is no
  $A$ such that $NS_{\kappa,\lambda}^{[\delta']^{<\theta'}} =
NS_{\kappa,\lambda}^{[\delta]^{<\theta}} | A$,  thus providing 
partial converses to Lemma 7.2.

\bigskip
\pmb{\bf Definition.} \ {\it  Assume $\overline\theta < \kappa$. Then for $f :
P_{\overline\theta\cdot 3} (\delta) \longrightarrow P_\kappa(\lambda)$ and
$X\subseteq\lambda$, we define $\Gamma_f(X)$ as follows. Let $\rho = 
\overline\theta \cdot
\aleph_0$ if $\overline\theta \cdot \aleph_0$ is a regular cardinal, 
and $\rho =
(\overline\theta \cdot \aleph_0)^+$ otherwise. Define $X_\alpha 
\subseteq \lambda$ for
$\alpha < \rho$ by~:
\begin{enumerate}
\item[$(0)$] \vskip-0,2cm $X_0 = X$.
\item[$(1)$]  \vskip-0,2cm $X_{\alpha+1} = X_\alpha \cup \big( \cup 
f[P_{\overline\theta\cdot 3}
(X_\alpha
\cap\delta)]\big)$.
\item[$(2)$]  \vskip-0,2cm $X_\alpha = \displaystyle\bigcup_{\beta < 
\alpha} X_\beta$ \ if
$\alpha$ is an infinite limit ordinal.
\end{enumerate}

Now let $\Gamma_f(X) = \displaystyle\bigcup_{\alpha <\rho} X_\alpha$.}

\vfill\eject
Notice that
\vskip-0,5cm
$$\Gamma_f(X) = \displaystyle\bigcap \ \{ Y : X\subseteq Y 
\subseteq\lambda \ \hbox{ and } \
(\forall e\in P_{\overline\theta\cdot 3} (Y\cap\delta)) \ f(e) \subseteq Y\}.$$

\bigskip
{\bf Definition.} \  {\it Let $\delta'$ be an ordinal with $\delta 
\leq\delta'\leq\lambda$,
and $\theta'$ be a cardinal with \break 
$\theta\leq\theta'\leq\kappa$. Given $f :
P_{\overline{\theta'}\cdot 3} (\delta') \longrightarrow 
P_\kappa(\lambda)$ and $k :
P_{\overline\theta\cdot 3} (\delta) \longrightarrow 
P_\kappa(\lambda)$, we define \break $u(f, k)
: P_{\overline{\theta'}\cdot 3} (\delta') \longrightarrow 
P_\kappa(\lambda)$ by~:
$(u(f,k)) (e) = f(e) \cup k(e)$ if $e \in P_{\overline\theta\cdot 3} 
(\delta)$, and
\break $(u(f,k)) (e) = f(e)$  otherwise.}

\smallskip
Notice that if $\overline{\theta'} < \kappa$ and there exists a 
$[\delta']^{<\theta'}$-normal
ideal on $P_\kappa(\lambda)$,  then \break $\Gamma_{u(f,k)} (a) \in 
C_f^{\kappa,\lambda} \cap
C_k^{\kappa,\lambda}$ for every $a\in P_\kappa (\lambda)$ with 
$\overline{\theta'}\cdot 3
\subseteq a$.\\

\bigskip
\hskip-0,2cm\pmb{\sc  Proposition 7.6.} \  {\it Let $\delta'$ be an 
ordinal with
$\kappa \cup\delta\leq\delta'\leq\lambda$, and $\theta'$ be a cardinal with
$\theta\leq\theta'\leq\kappa$. Assume that
$|\delta|^{<\overline\theta} < |\delta'|^{<\overline{\theta'}} < 
\lambda$.  Then
$NS_{\kappa,\lambda}^{[\delta']^{<\theta'}} \not= 
NS_{\kappa,\lambda}^{[\delta]^{<\theta}} |
A$ for all $A \in {\big( NS_{\kappa,\lambda}^{[\delta]^{<\theta}}\big)}^+$.
}\\

\vspace*{0ex}{\it Proof.} Let $f : P_{\overline{\theta'}\cdot 
3}(\delta')\longrightarrow
P_\kappa (\lambda)$. Set \  $\nu = \kappa\cup 
{(|\delta]^{<\overline\theta})}^+$ and select a
one-to-one

$i : \nu\longrightarrow P_{\overline{\theta'}\cdot 3} (\delta')$ and 
a one-to-one
$j : \nu\longrightarrow \lambda - (\nu\cup\delta\cup 
(\cup~\hbox{rang} (f)))$.  Define $h :
P_{\overline{\theta'}\cdot 3} (\delta') \longrightarrow P_2(\lambda)$ 
so that $h(i(\xi)) = \{
j(\xi)\}$ for every $\xi\in\nu$. Now let $k : P_{\overline\theta\cdot 
3} (\delta)
\longrightarrow P_\kappa(\lambda)$. Pick $\xi\in\nu$ so that 
$j(\xi)\notin\cup~ \hbox{ran}
(k)$.

First assume $\overline{\theta'}<\kappa$. We set $b = \Gamma_{u(f,k)} 
((\overline{\theta'}\cdot 3)
\cup i (\xi))$. Then $b\notin C_h^{\kappa,\lambda}$ since 
$j(\xi)\notin b$. Next assume
$\overline{\theta'} = \kappa$. We define $d_\beta\in P_\kappa(\lambda)$ and
$\gamma_\beta\in\kappa$ for $\beta < \kappa$ as follows~:

\begin{enumerate}
\item[(0)] \vskip-0,2cm $d_0 = \{ 0\} \cup i(\xi) \cup |i(\xi)|^+$ if 
$\overline\theta = \kappa$
and
$d_0 = (\overline\theta\cdot 3) \cup i(\xi) \cup |i(\xi)|^+$ otherwise.
\item[(1)] \vskip-0,2cm $\gamma_\beta = \cup (d_\beta \cap \kappa)$.
\item[(2)] \vskip-0,2cm $d_{\beta+1} = d_\beta \cup (\gamma_\beta+1) 
\cup (\cup \{ (u(f,k)) (e) :
e\in P_{|d_\beta\cap\kappa|} (d_\beta \cap \delta') \} )$.
\item[(3)] \vskip-0,2cm $d_\beta = \displaystyle\bigcup_{\zeta 
<\beta} d_\zeta$ if $\beta$ is an
infinite limit ordinal.
\end{enumerate}

Select a regular infinite cardinal $\tau <\kappa$ so that  (a) 
$\gamma_\tau = \tau$, and (b)
$\overline\theta\leq\tau$ if $\overline\theta <\kappa$.

Then $d_\tau\in C_f^{\kappa,\lambda} \cap C_k^{\kappa,\lambda}$. 
Moreover  $i(\xi) \in
P_{|d_\tau\cap(\overline{\theta'}\cdot 3)|}  (d_\tau\cap\delta')$ and 
$j(\xi)\notin d_\tau$, so
$d_\tau\notin C_h^{\kappa,\lambda}$.\eop

\bigskip
\hskip-0,2cm\pmb{\sc  Proposition 7.7.} \  {\it Let $\mu$ be a 
cardinal with $\kappa\leq\mu
<\lambda$. Assume that either $\lambda$ is a regular cardinal, or 
$u(\mu^+,\lambda) =
\lambda$. Then $NS_{\kappa,\lambda} \not= NS_{\kappa,\lambda}^\mu | 
A$ for every $A\in
{(NS_{\kappa,\lambda}^\mu)}^+$.
}\\
\vskip-0,4cm
\vspace*{0ex}{\it Proof. } Let us first deal with the case when 
$\lambda$ is regular. Fix $f
: P_3(\lambda) \longrightarrow P_\kappa(\lambda)$. Let $C$ be the set of all
$\beta\in\lambda$ such that  $f(e)\subseteq\beta$ for every $e\in 
P_3(\beta)$. Notice that
$C$ is a closed unbounded set. Define $h : P_2(\lambda) 
\longrightarrow P_2(\lambda)$ so
that $h(\{\xi\}) = \{ \beta_\xi\}$, where $\beta_\xi$ is the least 
element $\beta$ of $C$
such that $\beta > 3\cup\xi$. Now given $k : P_3(\mu) \longrightarrow 
P_\kappa(\lambda)$,
select $\xi\in\lambda$ so that $\cup ~\hbox{ran} (k) \subseteq\xi$. 
Setting $b =
\Gamma_{u(f,k)} (3\cup \{ \xi\})$, we have $b\notin 
C_h^{\kappa,\lambda}$ since $h(\{ \xi\})
\not\subseteq b$.

Next assume that $\lambda$ is a singular cardinal and $u(\mu^+, 
\lambda) = \lambda$. Fix $f
: P_3(\lambda) \longrightarrow P_\kappa(\lambda)$. Select a one-to-one $j :
\lambda\longrightarrow P_{\mu^+} (\lambda)$ so that  $\hbox{ran}~ (j)\in
I_{\mu^+,\lambda}$. Define  $h : P_2(\lambda) \longrightarrow 
P_2(\lambda)$ so that $h(\{
\xi\}) = \{ \beta_\xi\}$, where $\beta_\xi$ is the least element 
$\beta$ of $\lambda$ such
that $\beta\notin \Gamma_f(\{\xi\} \cup j(\xi))$. Now given $k : 
P_3(\mu) \longrightarrow
P_\kappa(\lambda)$, select $\xi\in\lambda$ so that $3\cup 
(\cup~ \hbox{ran} (k)) \subseteq
j(\xi)$. Set $b = \Gamma_{u(f,k)} (3\cup \{\xi\})$. Then $b\subseteq 
\Gamma_f (\{ \xi\} \cup
j(\xi))$ and therefore $b\notin C_h^{\kappa,\lambda}$. \eop

\hskip-0,2cm\pmb{\sc  Proposition 7.8.} \  {\it Let $\sigma$ be a 
cardinal such that {\rm (a)}
$\kappa < \sigma\leq\lambda$, and letting $\theta = (cf (\sigma))^+$, 
{\rm (b)} $\theta <\kappa$,
{\rm (c)} $\sigma^{<\theta}\geq\lambda$, {\rm (d)} $\mu^{<\theta} < 
\sigma$ for every cardinal $\mu <
\sigma$, and {\rm (e)} $u(\sigma,\lambda)\leq\lambda^{<\theta}$. 
Further let $\nu$ be a cardinal
with $\kappa\leq\nu <\sigma$. Then 
$NS_{\kappa,\lambda}^{[\sigma]^{<\theta}} \not=
NS_{\kappa,\lambda}^{[\nu]^{<\theta}} | A$  for every \break $A\in
{(NS_{\kappa,\lambda}^{[\nu]^{<\theta}})}^+$.
}
\bigskip

\vspace*{0ex}{\it Proof. } Fix $f : P_\theta (\sigma)\longrightarrow 
P_\kappa(\lambda)$.
Select $A\in I_{\sigma,\lambda}^+$ so that $A\subseteq \{ a\in 
P_\sigma(\lambda) :
\kappa\subseteq a\}$ and $|A| \leq \lambda^{<\theta}$. From Lemma 7.3 we get
$\lambda^{<\theta} = \sigma^{<\theta}$. So we can find a one-to-one 
\break $j : A\longrightarrow
P_\theta (\sigma)$. Notice that if $a\in A$, then setting $\mu = 
|a\cup j(a)|$, we have \break
$|\Gamma_f(a \cup j(a))| \leq \mu^{<\theta}$ since by Proposition 2.23
${(\mu^{<\theta})}^{<\theta} = \mu^{<\theta}$. Define $h : P_\theta (\sigma)
\longrightarrow P_2(\lambda)$ so that for every $a\in A$, $h(j(a)) = 
\{ \xi_a\}$, where
$\xi_a$ is the least element of the set $\lambda - \Gamma_f (a\cup 
j(a))$. Now given $k :
P_\theta (\nu)\longrightarrow P_\kappa(\lambda)$, pick $a\in A$ so 
that $\cup~\hbox{ran}
(k) \subseteq a$, and put $b = \Gamma_{u(f,k)} (\theta\cup j(a))$. Then
$h(j(a))\not\subseteq b$ since $b\subseteq \Gamma_f (a\cup j(a))$, 
hence $b\notin
C_h^{\kappa,\lambda}$. \eop

\hskip-0,2cm\pmb{\sc  Corollary  7.9.} \  {\it Assume that $\theta = 
(cf (\lambda))^+$,
$\theta < \kappa$ and $\mu^{<\theta} < \lambda$ for every cardinal 
$\mu < \lambda$. Then
for every cardinal $\nu$ with $\kappa\leq\nu <\lambda$, and every $A\in
{\big(NS_{\kappa,\lambda}^{[\nu]^{<\theta}}\big)}^+$,\break
$NS_{\kappa,\lambda}^{[\lambda]^{<\theta}} \not= 
NS_{\kappa,\lambda}^{[\nu]^{<\theta}} | A$.
}

\eject
{\Large\bf 8. Dominating numbers}\\

{\bf Throughout this section $\mu$ will denote a cardinal with $\mu >0$}.

\bigskip
The dominating numbers we will consider now are three-dimensional 
generalizations of the well-known
cardinal invariant $\dom$. The connection with the notion of 
$\dt$-normality will be established in the
next section.

\bigskip
We will see that our numbers admit several equivalent definitions. It 
is convenient to give the
following `unofficial' definition first.\\

{\bf Definition.} \  {\it $\deklmu$ is the smallest cardinality of 
any $F\subseteq\emupkl$
such that  for every $g\in\emupkl$, there is $f\in F$ with
$|\{\alpha\in\mu:g(\alpha)\not\subseteq f(\alpha)\}| <\mu$.}

\bigskip
The following two propositions will be very useful.\\

\bigskip
\hskip-0,2cm\pmb{\sc Proposition 8.1.} \  {\it  $\deklmu\ge 
u(\kappa,\lambda)$.}\\
\vspace*{2ex}{\it Proof. }Given $F\subseteq\emupkl$ with $|F|<
u(\kappa,\lambda)$, it is easy to define
$g\in\emupkl$ so  that $g(\alpha)\not\subseteq f(\alpha)$ for all 
$\alpha\in\mu$ and $f\in F$.\eop

\hskip-0,2cm\pmb{\sc Proposition 8.2.} {\it $cf(\deklmu)>\mu$.}\\
\vspace*{0ex}{\it Proof. }We can assume that $\mu\ge\omega$, since the
result is immediate from Propositions 8.1 and 1.3 (0)  if
$\mu<\omega$. Select a bijection $j:\mu\times\mu\longrightarrow\mu$.
Suppose toward a
\mbox{contradiction} that there are
$F_\gamma\subseteq\emupkl$ for $\gamma<\mu$ such that a) 
$|F_\gamma|<\deklmu$ for all
$\gamma<\mu$, b) $F_\gamma\cap  F_\delta=\emptyset$ whenever 
$\gamma,\delta<\mu$ are such
that $\gamma\not=\delta$, and c) for every $g\in\emupkl$,  there is
$f\in\displaystyle\bigcup_{\gamma<\mu}F_\gamma$ with
$|\{\alpha<\mu:g(\alpha)\not\subseteq f(\alpha)\}|<\mu$. For  each
$\gamma<\mu$, there is $g_\gamma\in\emupkl$ such that
$|\{\alpha<\mu:g_\gamma(\alpha)\not\subseteq f(j(\gamma,\alpha))\}|=\mu$
  for every $f\in F_\gamma$. Define $h\in\emupkl$ by
$h(j(\gamma,\alpha))=g_\gamma(\alpha)$. There are
$\gamma<\mu$ and $f\in F_\gamma$ such that 
$|\{\beta<\mu:h(\beta)\not\subseteq f(\beta)\}|<\mu$. Then 
$|\{\alpha<\mu:h(j(\gamma,\alpha))\not\subseteq
  f(j(\gamma,\alpha))\}|<\mu$, a contradiction.\eop

{\bf Definition.} \ {\it $F\subseteq\emupkl$ is {\it 
$\emupkl$-dominating} if for every
$g\in\emupkl$, there is
$f\in F$ such that $g(\alpha)\subseteq f(\alpha )$ for all $\alpha<\mu$.}

\bigskip
The `official' definition of our three-cardinal version of the
dominating number $\dom$ reads as follows.

\bigskip
{\bf Definition.} \ {\it $\dklmu$ is the least cardinality of any
$\emupkl$-dominating
$F\subseteq\emupkl$.}

\smallskip
Let us first observe that $\dklmu$ is a familiar quantity in case
$\mu<\kappa$.\\

\hskip-0,2cm\pmb{\sc Proposition 8.3.} \  {\it Assume $\mu<\kappa$. Then
$\dklmu=  u(\kappa,\lambda)$.}

\smallskip
\vspace*{0ex}{\it Proof. }Since clearly $\dklmu\ge\deklmu$, we get
$\dklmu\ge u(\kappa,\lambda)$ by Proposition 8.1.  For the reverse
inequality, observe  that given $g\in\emupkl$, we have
$g(\alpha)\subseteq\cup \ {\rm ran}(g)$ for all $\alpha<\mu$.\eop

\bigskip
\hskip-0,2cm\pmb{\sc Proposition 8.4.} \  {\it $\dklmu=\deklmu$.}

\smallskip
\vspace*{0ex}{\it Proof. }It is immediate that $\dklmu\ge\deklmu$.
If $\mu<\kappa$, the reverse inequality follows from
Propositions 8.1 and 8.3.

\bigskip
Now assume  $\mu\ge\kappa$. Select a bijection
$j:\mu\times\mu\longrightarrow\mu$, and
let $F\subseteq\emupkl$ be such that for every $g\in\emupkl$,
there is $f\in F$ with $|\{\alpha<\mu:g(\alpha)\not\subseteq
f(\alpha)\}|<~\mu$.\break  For $f\in F$ and $\beta<\mu$,
define $f_\beta\in\emupkl$ by  $f_\beta(\xi)=f(j(\beta,\xi))$. Notice
that \break
$|\{f_\beta:\beta<\mu$ and $f\in F\}|\le|F|$ by Proposition
8.2. Given $h\in\emupkl$, define $g\in\emupkl$ by
$g(j(\beta,\xi))=h(\xi)$. Pick $f\in F$ with
$|\{\alpha<\mu:g(\alpha)\not\subseteq f(\alpha)\}|<\mu$.
There exists $\beta<\mu$ such that
\vskip-0,2cm
\centerline{$\{\alpha<\mu:g(\alpha)\not\subseteq 
f(\alpha)\}\cap\{j(\beta,\xi):\xi<\mu\}=\emptyset$.}
  \vskip-0,3cm
Then
\vskip-0,3cm
\centerline{$h(\xi)=g(j(\beta,\xi))\subseteq f(j(\beta,\xi))=f_\beta(\xi)$}

  for every $\xi<\mu$.\eop

The following will be repeatedly used.\\

\vskip-0,5cm
\hskip-0,2cm\pmb{\sc Corollary 8.5.} \  {\it $\dklmu\ge\lambda$.} \\
\vspace*{2ex}{\it Proof. }By Propositions 8.4, 8.1 and 1.3 (0).\eop

\vskip-1cm
Let us consider another variation on the definition of $\dklmu$.\\

{\bf Definition.} \ {\it $\Delta_{\kappa,\lambda}^\mu$ is the least 
cardinality of any
$F\subseteq\emupkl$ with the property that for every
$g\in{^\mu}{\lambda}$, there is $f\in F$ such that $g(\alpha)\in 
f(\alpha)$ for all $\alpha\in\mu$.}
\\

\hskip-0,2cm\pmb{\sc Proposition 8.6.} \  {\it
$\Delta_{\kappa,\lambda}^\mu\le\dklmu\le\Delta_{\kappa,\lambda}^{\mu\cdot\tau}$, 
where $\tau=\kappa$  if $\kappa$ is a limit cardinal, and $\tau=\nu$ if
$\kappa=\nu^+$.} \\
\vspace*{0ex}{\it Proof. }It is immediate that
$\Delta_{\kappa,\lambda}^\mu\le\dklmu$. Let us show the other
inequality. Select a bijection $j_a:|a|\longrightarrow
  a$ for each $a\in\pkl$. Let $F\subseteq{^{(\mu\times\tau)}}\!{\pkl}$ 
be such that for every $g\in
{^{(\mu\times\tau)}}\!{\lambda}$, there is $f\in F$ with the property 
that $g(\gamma,\xi )\in
f(\gamma,\xi )$  for every $(\gamma,\xi )\in\mu\times\tau$. For 
$f\in F$, define $k_f\in\emupkl$
by
$k_f(\gamma)=\displaystyle\bigcup\{f(\gamma,1+\xi ):\xi 
<\cup(\kappa\cap f(\gamma,0))\}$.
Given $h\in\emupkl$, define $g\in{^{(\mu\times\tau)}}\!{\lambda}$ by~: (i)
$g(\gamma,0)=|h(\gamma)|$, and
(ii) $g(\gamma,1+\xi)=j_{h(\gamma)}(\xi)$ if $\xi<g(\gamma,0)$, and
$g(\gamma,1+\xi)=0$ otherwise.  There is $f\in F$ such that
$g(\gamma,\xi)\in f(\gamma,\xi)$ for all
$(\gamma,\xi)\in\mu\times\tau$.  We have that $h(\gamma)\subseteq
k_f(\gamma)$ for every $\gamma\in\mu$. Hence $\dklmu\le|F|$. \eop

We will now see that $\dklmu$ is easy to compute if $\lambda$ is large with
respect to $\mu$.\\

\hskip-0,2cm\pmb{\sc Lemma 8.7.} \
\begin{enumerate}
\item[(0)] \vskip-0,3cm {\it Assume $\mu<\kappa$. Then
$\lambda^{<\kappa}=\dklmu\cdot 2^{<\kappa}$.}
\item[(1)] \vskip-0,3cm {\it Assume $\mu\ge\kappa$. Then $\lambda^\mu=\dklmu\cdot 2^\mu$.}
\end{enumerate}

\vspace*{0ex}{\it Proof. }\vskip-0,2cm
\begin{enumerate}
\vskip-0,3cm
\item[(0)] : It is well-known (see [DoM]) that
$\lambda^{<\kappa}=u(\kappa,\lambda)\cdot 2^{<\kappa}$. So the result 
follows from
Proposition 8.3.
\item[(1)] : $\lambda^\mu=|\emupkl|\le\dklmu\cdot 
|{^{\mu}}{(2^{<\kappa})}|\le|\emupkl|$.\eop
\end{enumerate}

\hskip-0,2cm\pmb{\sc Proposition 8.8.} \
\begin{enumerate}
\item[(0)]\vskip-0,2cm {\it Assume  $\mu<\kappa$ and $\lambda\ge 
2^{<\kappa}$. Then
$\dklmu=\lambda^{<\kappa}$.}
\item[(1)]\vskip-0,2cm {\it Assume  $\mu\ge\kappa$ and $\lambda\ge 2^\mu$.
Then $\dklmu=\lambda^\mu$.}
\end{enumerate}
\vspace*{0ex}{\it Proof. }By Lemma 8.7 and Corollary 8.5.\eop

\hskip-0,2cm\pmb{\sc Proposition 8.9.} \  {\it Assume}  GCH. {\it Then}
\begin{enumerate}
\item[a)] $\dklmu=\mu^+$ if $\mu\ge\lambda$.

\item[b)] $\dklmu=\lambda^+$ if $\mu<\lambda$ and $\mu^+\cdot 
\kappa>cf(\lambda)$.

\item[c)] $\dklmu=\lambda$ if $\mu^+\cdot \kappa\le cf(\lambda)$.
\end{enumerate}
\vspace*{0ex}{\it Proof. }a) : By Propositions 8.2 and 8.4 and Lemma 8.7 (1).

\hskip1,2cm b) and c) : By Proposition 8.8.\eop

Notice that $\dklmu\ge\lambda$ and $cf(\dklmu)\ge\mu^+\cdot \kappa$ 
by Corollary
8.5 and Propositions 8.2, 8.4,  8.3 and 1.3 (1). Thus Proposition 8.9 
shows that
$\dklmu$ assumes its least possible value under
GCH. Let us now show that $\kappa$-c.c. forcing preserves this 
minimal value in case
$\kappa>\omega$.\\

\bigskip
\hskip-0,2cm\pmb{\sc Proposition 8.10.} \  {\it Assume 
$\kappa>\omega$, and let $(P,<)$ be
a
$\kappa$-c.c. notion of forcing. Then
$(\dom_{\kappa,\lambda}^{|\mu|})^{V^P}\le(\dklmu)^V$.}\\
\vspace*{0ex}{\it Proof. }Let $G$ be $P$-generic over $V$. Given an 
ordinal $\xi$ and
$f:\xi\longrightarrow\lambda$ in $V\lbrack  G\rbrack$, there
is by Lemma 6.8 in chapter VII of [K], $F:\xi\longrightarrow\pkl$ in 
$V$ with the property that
$f(\alpha)
\in F(\alpha)$ for every $\alpha<\xi$. It immediately follows that 
$(\Delta_{\kappa,\lambda}^{|\mu|})^{
V\lbrack G\rbrack}\le(\dklmu)^V$, which by Proposition~8.6 gives
$(\dom_{\kappa,\lambda}^{|\mu|})^{V\lbrack G\rbrack}\le(\dklmu)^V$ if
$\mu\ge\kappa$.

\bigskip
Now assume  $\mu<\kappa$. Then $(\dom_{\kappa,\lambda}^{|\mu|})^{V\lbrack
G\rbrack}=(u(\kappa,\lambda))^{ V\lbrack G\rbrack}$ and 
$(\dklmu)^V=(u(\kappa,\lambda))^V$
by Proposition 8.3. In $V$, let $A\in\iklp$. In $V\lbrack  G\rbrack$, 
let $b\in\pkl$, and
select a bijection $j:|b|\longrightarrow b$. There exists 
$F:|b|\longrightarrow\pkl$ in
$V$ such that $j(\alpha)\in F(\alpha)$ for all $\alpha<|b|$. Pick 
$a\in A$ with $\cup$ ran$(F)\subseteq a$.
Then $b\subseteq a$. Thus it is still true in $V\lbrack G\rbrack$ 
that $A\in\iklp$. It follows that $(u(\kappa,\lambda)
)^{V\lbrack G\rbrack}\le(u(\kappa,\lambda))^V$.\eop

We will present a few identities and inequalities that can be used to evaluate
$\dklmu$ in the absence of GCH.
The following is immediate.\\

\bigskip
\hskip-0,2cm\pmb{\sc Lemma 8.11.} \  {\it Let $\tau$ and $\nu$ be 
cardinals such that
$\tau\ge\lambda$ and $\nu\ge\mu$. Then
$\dom_{\kappa,\tau}^\nu\ge\dklmu$.}\\

\bigskip
\hskip-0,2cm\pmb{\sc Proposition 8.12.} \  {\it Assume  $\lambda>\kappa$ and
$cf(\lambda)\ge\kappa\cdot\mu^+$. Then
$\dklmu=\lambda\cdot
(\displaystyle\bigcup_{\kappa\le\rho<\lambda}
\dom_{\kappa,\rho}^\mu)$.}\\
\vspace*{2ex}{\it Proof. }$\le$ : Observe that $\emupkl=
\displaystyle\bigcup_{\kappa\le\alpha<\lambda}\emupka$.

\bigskip
\hskip1,3cm $\ge$ : By Corollary 8.5 and Lemma 8.11.\eop

We will use the following two-cardinal version of $\dom$.\\

{\bf Definition.} \  {\it $\dom_\lambda^\mu$ is the least cardinality of any
$F\subseteq{^\mu}{\lambda}$ with the  property that for every $g\in 
{^\mu}{\lambda}$, there
is $f\in F$ such that
$g(\alpha)\le f(\alpha)$ for every $\alpha<\mu$.

We put $\dom_\kappa=\dom_\kappa^\kappa$.}

\bigskip
Thus $\dom=\dom_\omega$.\\

\bigskip
\hskip-0,2cm\pmb{\sc Lemma 8.13.} \  {\it Assume $cf(\lambda)\ge\kappa$. Then
$\Delta_{\kappa,\lambda}^\mu\ge\dom_\lambda^\mu$.}\\
\vspace*{0ex}{\it Proof. }Let $F\subseteq\emupkl$ be such that for
every $g\in{^\mu}{\lambda}$, there is $f\in F$ with the  property
that $g(\alpha)\in f(\alpha)$ for all $\alpha<\mu$. Given
$g\in{^\mu}{\lambda}$, select $f\in F$ so that
$g(\alpha)\in f(\alpha)$ for all $\alpha<\mu$. Then $g(\alpha)\le\cup 
f(\alpha)$ for every
$\alpha<\mu$.\eop

\hskip-0,2cm\pmb{\sc Proposition 8.14.} \  {\it 
$\dom_{\kappa,\kappa}^\mu=\dom_\kappa^\mu$.}\\

\vspace*{0ex}{\it Proof. }We have
$\dom_{\kappa,\kappa}^\mu\ge\dom_\kappa^\mu$ by Lemma 8.13. Now let
$F\subseteq{^\mu}{\kappa}$ be such that
  for every $g\in{^\mu}{\kappa}$, there is $f\in F$ with the property 
that $g(\alpha)\le f(\alpha)$
for every $\alpha\in\mu$. Given $h\in\emupkk$, select $f\in F$ so that $\cup
h(\alpha)<f(\alpha)$  for all $\alpha\in\mu$. Then $h(\alpha)\subseteq
f(\alpha)$ for every $\alpha\in\mu$. Hence
$\dom_{\kappa,\kappa}^\mu\le
\dom_\kappa^\mu$.\eop

The following basic observation is very fruitful.\\

\hskip-0,2cm\pmb{\sc Proposition 8.15.} \
\begin{enumerate}
\item[(0)] $\dklmu\le
\dom_{\kappa,\rho}^\mu\cdot\dom_{\rho^+,\lambda}^\mu\le\dom_{\kappa,\lambda}^{\mu
\cdot\rho}$
{\it for every cardinal}
$\rho$ with $\kappa\le\rho<\lambda$.

\item[(1)] $\dklmu\le\dom_{\kappa,\rho}^\mu\cdot\dom_{\rho,\lambda}^\mu\le
\dom_{\kappa,\lambda}^{\mu\cdot\rho}$
{\it for  every regular cardinal}
$\rho$ with $\kappa\le\rho\le\lambda$.
\end{enumerate}
\vspace*{0ex}{\it Proof. }Fix a cardinal $\rho$ with
$\kappa\le\rho\le\lambda$, and let $\tau$ be a regular cardinal with
$\rho\le\tau\le\lambda\cap\rho^+$. Pick a bijection 
$j_a:|a|\longrightarrow a$ for each $a\in P_\tau(\lambda)$.

\bigskip
Let us first show that
$\dklmu\le\dom_{\kappa,\rho}^\mu\cdot\dom_{\tau,\lambda}^\mu$. Select a
${^\mu}\!{P_\kappa(
\rho)}$-dominating $F\subseteq{^\mu}\!{P_\kappa(\rho)}$ and a 
${^\mu}\!{P_\tau(\lambda)}$-dominating
$G\subseteq{^\mu}\!{P_\tau(\lambda)}$.
Define $\varphi:F\times G\longrightarrow\emupkl$
by

\centerline{$(\varphi(f,g))(\alpha)=j_{g(\alpha)}\lbrack 
f(\alpha)\cap|g(\alpha)|\rbrack$.}

  We claim that ran$(\varphi)$ is $\emupkl$-dominating. Thus let 
$r\in\emupkl$. Pick $g\in G$ so that
$r(\alpha)\subseteq g(\alpha)$ for all $\alpha<\mu$. Then pick $f\in 
F$ so that $j_{g(\alpha)}^{-1}(r(\alpha))\subseteq
  f(\alpha)$ for every $\alpha<\mu$. Then 
$r(\alpha)\subseteq(\varphi(f,g))(\alpha)$ for all $\alpha<\mu$, 
which proves our
  claim.

\bigskip
Let us next show that
${\dom}_{\kappa,\rho}^\mu\cdot\dom_{\tau,\lambda}^\mu\le
\dom_{\kappa,\lambda}^{\mu\cdot\rho}$.
We have
  $\dom_{\kappa,\rho}^\mu\le\dom_{\kappa,\lambda}^{\mu\cdot\rho}$ by 
Lemma 8.11. Now let
$H\subseteq   {^{(\mu\times\rho)}}\!{\pkl}$ be such that for every
$p\in{^{(\mu\times\rho)}}\!{\pkl}$, there is
$h\in H$ with the property that $p(\alpha,\beta)\subseteq 
h(\alpha,\beta)$ for every $(\alpha,\beta)\in\mu\times\rho$.
Given $q\in {^\mu}\!{P_\tau(\lambda)}$, pick $h\in H$ so that 
$\{j_{q(\alpha)}(\beta)\}\subseteq h(\alpha,\beta)$
whenever $\alpha\in\mu$ and $\beta\in|q(\alpha)|$. If $\tau=\rho^+$, then
$q(\alpha)\subseteq\displaystyle\bigcup_{\beta\in\rho}h(\alpha,\beta)
$, and we can conclude that
$\dom_{\tau,\lambda}^\mu\le\dom_{\kappa,\lambda}^{\mu\cdot\rho}$. Now assume
$\tau=\rho$, and let $K\subseteq {^\mu}{\tau}$ be such that for every
$i\in{^\mu}{\tau}$,  there is $k\in K$ with the property that $i(\alpha)\le
k(\alpha)$ for all $\alpha<\mu$. Then there is $k\in K$ such  that 
$|q(\alpha)|\le
k(\alpha)$ for every $\alpha<\mu$. We have
$q(\alpha)\subseteq\displaystyle\bigcup_{\beta\in
  k(\alpha)}h(\alpha,\beta)$ for all $\alpha<\mu$. Thus
$\dom_{\tau,\lambda}^\mu\le\dom_{\kappa,\lambda}^{\mu\cdot\rho}\cdot\dom_\tau^\mu
$, which gives 
$\dom_{\tau,\lambda}^{\mu}\le\dom_{\kappa,\lambda}^{\mu\cdot\rho}$, 
since
$\dom_\tau^\mu\le\dom_{\kappa,\tau}^\mu\le
\dom_{\kappa,\lambda}^{\mu\cdot\rho}$ by Lemmas 8.11 and 8.13 and 
Proposition 8.6.\eop

\hskip-0,2cm\pmb{\sc Corollary 8.16.} \  {\it Let $\sigma$ and $\chi$ 
be uncountable cardinals
such that
$\sigma\le\mu\cap\chi$ and
$cf(\sigma)=\omega$. Then there is a regular infinite cardinal 
$\tau<\sigma$ such
that $\dom_{\rho,\chi}^\mu =\dom_{\tau,\chi}^\mu$ for every regular 
cardinal $\rho$
with $\tau\le\rho<\sigma$.}\\
\vspace*{0ex}{\it Proof. }Pick regular infinite cardinals
$\sigma_0<\sigma_1<\sigma_2<\ldots$ with
$\sigma=
\displaystyle\bigcup_{n\in\omega}\sigma_n$. Then by Proposition 8.15 (1),
$\dom_{\sigma_0,\chi}^\mu\ge\dom_{\sigma_1,\chi}^\mu\ge\dom_{\sigma_2, 
\chi}^\mu\ge\ldots\
\!$,   ~and
$\dom_{\sigma_n,\chi}^\mu\ge\dom_{\rho,\chi}^\mu\ge\dom_{\sigma_{n+1}, 
\chi}^\mu$
whenever $n\in\omega$ and $\rho$ is a regular cardinal with 
$\sigma_n\le\rho\le\sigma_{n+1}$. Hence there exists $
k\in\omega$ such that $\dom_{\sigma_k,\chi}^\mu=\dom_{\rho,\chi}^\mu$ 
for every regular cardinal $\rho
$ with $\sigma_k\le\rho<\sigma$.\eop

\hskip-0,2cm\pmb{\sc Corollary 8.17.} \  {\it Assume 
$\kappa\le\mu<\lambda$. Then
$\dklmu=\dom_{\kappa,\mu}^\mu\cdot u(\mu^+,\lambda)$.}\\
\vspace*{2ex}{\it Proof. }By Propositions 8.15 (0) and 8.3.\eop
\vskip-0,5cm
Proposition 8.3 and Corollary 8.17 show that for $\mu\le\lambda$, the value of
$\dklmu$ is determined by the values
  taken by $\dom_{\kappa,\tau}^\tau$ and $u(\tau,\lambda)$ when $\tau$ 
ranges from $\kappa$ to $\lambda$.

\bigskip
Let us next consider the relationship between $\dklmu$ and 
$\dom_{\kappa,\lambda^+}^\mu$.\\

\hskip-0,2cm\pmb{\sc Proposition 8.18.} \

\begin{enumerate}
\item[(0)] \vskip-0,3cm
$\dom_{\kappa,\lambda^{+n}}^\mu=\dklmu\cdot 
(\prod_{i=1}^n{\dom}_{\lambda^{+i}}^\mu)$
{\it for each} $n\in\omega-\{0\}$.

\item[(1)] \vskip-0,2cm {\it Assume}  $\mu\le\lambda$.  {\it Then} \
$\dom_{\kappa,\lambda^{+n}}^\mu=\dklmu\cdot \lambda^{+n}$
{\it for every} $n\in\omega$.
\end{enumerate}
\vspace*{0ex}{\it Proof. }
\begin{enumerate}
\item[(0)] : We get
$\dom_{\kappa,\lambda^+}^\mu=\dklmu\cdot\dom_{\lambda^+}^\mu$ by 
Propositions 8.15 (0),  8.14
and 8.6 and Lemmas 8.13 and 8.11. The desired result is then obtained 
by induction.
\item[(1)] : The result follows from (0) and Propositions 8.3 and 
8.14 if $n>0$, and from
Corollary 8.5 otherwise.\eop
\end{enumerate}

\hskip-0,2cm\pmb{\sc Corollary 8.19.} \
\begin{enumerate}
\item[(0)]\vskip-0,2cm $\dom_{\kappa,\kappa^{+n}}^\mu=\prod_{i=0}^{n} 
{\dom}_{\kappa^{+i}}^\mu$
{\it for every} \  $n\in\omega$.
\item[(1)]\vskip-0,2cm 
$\dom_{\kappa,\kappa^{+n}}^\kappa=\dom_\kappa\cdot\kappa^{+n}$ {\it 
for
every} \  $n\in\omega$.
\item[(2)] \vskip-0,2cm $\dom_{\kappa,\lambda^+}^\lambda=\dom_{\kappa,
\lambda}^\lambda$.
\end{enumerate}

\vspace*{0ex}{\it Proof. }
\vskip-0,5cm
\begin{enumerate}\vskip-0,5cm
\item[(0)] : By Propositions 8.18 (0) and 8.14.
\item[(1)] : By Propositions 8.18 (1) and 8.14.
\item[(2)] : By Propositions 8.18 (1), 8.4 and 8.2.\eop
\end{enumerate}

\vskip-0,5cm
\hskip-0,2cm\pmb{\sc Lemma 8.20.} \  {\it 
$\dklmu\le\dom_{\omega,u(\kappa,\lambda)}^\mu$.} \\

\vskip-0,5cm
\vspace*{0ex}{\it Proof. }If $\mu\ge\kappa$, then
$\dom_{\omega,u(\kappa,\lambda)}^\mu\ge\dom_{\kappa,u(\kappa,\lambda)} 
^\mu\ge\dklmu$ by
  Propositions 8.15 (1) and 1.3 (0) and Lemma 8.11. If $\mu<\kappa$, then
$\dom_{\omega,u(\kappa,
\lambda)}^\mu\ge u(\kappa,\lambda)=\dklmu$ by Corollary 8.5 and
  \linebreak Proposition 8.3.\eop
\vskip-0,5cm
\hskip-0,2cm\pmb{\sc Proposition 8.21.} \  {\it Assume 
$u(\kappa,\lambda)=\lambda$. Then
$\dom_{\omega,\lambda}^\mu=\dom_{\omega,\kappa}^\mu\cdot\dklmu$.}

\bigskip
\vspace*{0ex}{\it Proof. }By Lemmas 8.20 and 8.11 and Proposition
8.15 (1).\eop

Notice that 
$\dom_{\omega,\lambda}^1=\dom_{\omega,\kappa}^1\cdot\dom_{\kappa,\lambda}^1$
 if
and only if $u(\kappa,\lambda)=\lambda$.

Let us now deal with the computation of $\dom_{\kappa,\lambda^{<\eta}}^\mu$.\\
\vskip-0,5cm
\hskip-0,2cm\pmb{\sc Proposition 8.22.} \vskip-0,5cm
\begin{enumerate}\vskip-0,5cm
\item[(0)] 
$\dom_{\kappa,\lambda^{<\eta}}^\mu=\dom_{\kappa,\lambda\cdot\kappa^{<\ 
eta}}^\mu$
{\it  for every cardinal} $\eta$  {\it with} $\omega<\eta<\kappa$.
\item[(1)]\vskip-0,3cm
$\dom_{\kappa,\lambda^{<\eta}}^\mu=\dom_{\kappa,2^{<\eta}}^\mu\cdot\do 
m_{\eta,\lambda}^\mu$
{\it for every regular cardinal}
$\eta$ {\it with} $\kappa\le\eta\le\lambda$.

\item[(2)] \vskip-0,3cm
$\dom_{\kappa,\lambda^{<\eta}}^\mu=\dom_{\kappa,\eta^{<\eta}}^\mu\cdot 
\dom_{\eta^+,\lambda}^\mu$
{\it for every singular cardinal}
$\eta$ {\it such that} $\kappa\le\eta<\lambda$ {\it and either} 
$\eta<\mu$, {\it or}
$\eta^+=\lambda$.

\item[(3)] \vskip-0,3cm
$\dom_{\kappa,\lambda^{<\eta}}^\mu=\dom_{\kappa,2^{<\eta}}^\mu\cdot\do 
m_{\eta^+,\lambda}^\mu$
{\it for every regular  cardinal} $\eta$ {\it such that} 
$\kappa\le\eta<\lambda$ {\it and
either}
$\eta\le\mu$, {\it} or $\eta^+=\lambda$.
\end{enumerate}

\vspace*{0ex}{\it Proof. }(0), (1) and (2) : Let $\eta$ be an
uncountable cardinal
$\le\lambda$. Let us assume that $\eta<\lambda$
if $\kappa\le\eta$ and $\eta$ is singular. We define $\rho$ and $\tau$ by :
\begin{enumerate}
\item[(i)] \vskip-0,3cm $\rho=\kappa$ and $\tau=\kappa^{<\eta}$ if 
$\eta<\kappa$.

\item[(ii)]\vskip-0,3cm $\rho=\eta$ and $\tau=2^{<\eta}$ if 
$\kappa\le\eta$ and $\eta$ is regular.

\item[(iii)] \vskip-0,3cm $\rho=\eta^+$ and $\tau=\eta^{<\eta}$ if 
$\kappa\le\eta$ and $\eta$ is
singular.
\end{enumerate}
Let $F\subseteq {^{\mu}}\!{P_\rho(\lambda)}$ be 
${^{\mu}}\!{P_\rho(\lambda)}$-dominating,
and
$K\subseteq{^{\mu}}\!{P_\kappa(\tau)}$ be 
${^{\mu}}\!{P_\kappa(\tau)}$-dominating. Fix\break
  a bijection $j:\lambda^{<\eta}
\longrightarrow P_\eta(\lambda)$. For  $f\in F$ and $\alpha\in\mu$, 
select a one-to-one \break
$i_{f,\alpha}: j^{-1}(P_\eta(f(\alpha)))\longrightarrow\tau$. Given
$h\in{^{\mu}}\!{P_\kappa(\lambda^{<\eta})}$,  pick $f\in F$ so that
$\cup (j\lbrack h(\alpha)\rbrack)\subseteq f(\alpha)$ for every 
$\alpha\in\mu$. Then pick $k\in K$ so that
$i_{f,\alpha}\lbrack h(\alpha)\rbrack\subseteq k(\alpha)$ for all 
$\alpha\in\mu$. Then $h(\alpha)\subseteq i_{f,\alpha}^{-1}(k(\alpha))$
  for every $\alpha\in\mu$. Hence
$\dom_{\kappa,\lambda^{<\eta}}^\mu\le\dom_{\kappa,\tau}^\mu\cdot
\dom_{\rho,\lambda}^\mu$.

\bigskip
Since $\tau\le\lambda^{<\eta}$, we have
$\dom_{\kappa,\lambda^{<\eta}}^\mu\ge\dom_{\kappa,\tau}^\mu$  by Lemma 8.11. If
$\rho=\kappa$, 
$\dom_{\kappa,\lambda^{<\eta}}^\mu\ge\dom_{\rho,\lambda}^\mu$ by 
Lemma~
8.11.
If $\rho=\lambda$,
$\dom_{\kappa,\lambda^{<\eta}}^\mu\ge\dklmu\ge\dom_\lambda^\mu=\dom_{\rho,
\lambda}^\mu$  by
Lemmas 8.11 and 8.13 and Propositions 8.6 and 8.14. If 
$\kappa<\rho<\lambda\cap\mu^+$,
$\dom_{\kappa,\lambda^{<\eta}}^\mu=\dom_{\kappa,\lambda^{<\eta}}^{\mu\cdot
\rho}\ge\dom_{\rho,\lambda^{<\eta}}^\mu\ge
\dom_{\rho,\lambda}^\mu$
by Proposition 8.15 (1) and Lemma 8.11.
Finally, if $\rho=\eta$ and $\rho>\mu$,
$\dom_{\kappa,\lambda^{<\eta}}^\mu\ge\lambda^{<\eta}\ge 
u(\rho,\lambda)=\dom_{\rho,\lambda}^\mu$
by Corollary 8.5 and Proposition 8.3. Thus if it is not the case that
$\mu\cdot\kappa<\rho=\eta^+<\lambda$,  then
$\dom_{\kappa,\lambda^{<\eta}}^\mu\ge\dom_{\kappa,\tau}^\mu\cdot\dom_{ 
\rho,\lambda}^\mu$.

\bigskip
\hskip0,5cm (3) : Let $\eta$ be a regular cardinal with 
$\kappa\leq\eta <\lambda$. Assume
$\eta\le\mu$.  We get
$\dom_{\kappa,\lambda^{<\eta}}^\mu=\dom_{\kappa,2^{<\eta}}^\mu\cdot
\dom_{\eta,\lambda}^\mu$ by
(1).  Moreover, 
$\dom_{\eta,\lambda}^\mu=\dom_{\eta,\eta}^\mu\cdot\dom_{\eta^+,\lambda 
}^\mu$
by Proposition 8.15 (0), and
$\dom_{\kappa,2^{<\eta}}^\mu\ge\dom_{\kappa,\eta}^\mu\ge\Delta_{\kappa 
,\eta}^\mu\ge
\dom_\eta^\mu= \dom_{\eta,\eta}^\mu$
by Lemmas 8.11 and  8.13  and Propositions 8.6 and 8.14. It follows that
$\dom_{\kappa,\lambda^{
<\eta}}^\mu=\dom_{\kappa,2^{<\eta}}^\mu\cdot\dom_{\eta^+,\lambda}^\mu$.

\bigskip
Assume now  $\eta^+=\lambda$ and $\eta>\mu$. Since
$(\eta^+)^{<\eta}=\eta^{<\eta}\cdot\eta^+$ and
$\eta^{<\eta}=2^{<\eta}$, we have by Lemma 8.11 that
$\dom_{\kappa,\lambda^{<\eta}}^\mu=
\dom_{\kappa,2^{<\eta}\cdot\eta^+}^\mu=\dom_{\kappa,2^{<\eta}}^\mu\cdot
\dom_{\kappa,\eta^+}^\mu$.
The desired conclusion now follows from the following three observations~:
a) $\dom_{\kappa,\eta^+}^\mu=\dom_{\kappa,\eta}^\mu\cdot\eta^+$ by 
Proposition 8.18 (1).
\break b) 
$\dom_{\kappa,2^{<\eta}}^\mu\cdot\dom_{\kappa,\eta}^\mu=\dom_{\kappa,2 
^{<\eta}}^\mu$
by Lemma 8.11.
c) $\eta^+=\dom_{\eta^+,\lambda}^\mu$ by Propositions 8.3 and 1.3 (0).\eop

Let us make the following remark concerning Proposition 8.22 (3). 
Assume  GCH, and let
$\eta$ be a cardinal such that 
$\kappa\cdot\mu^+\le\eta=cf(\lambda)<\lambda$. Then
$\dom_{\kappa,\lambda^{<\eta}}^\mu\not=\dom_{\kappa,2^{<\eta}}^\mu\cdot
\dom_{\eta^+,\lambda}^\mu$,
since by Proposition 8.9
$\dom_{\kappa,\lambda^{<\eta}}^\mu=\lambda$ and 
$\dom_{\eta^+,\lambda}^\mu=\lambda^+$.\\

\bigskip
\hskip-0,2cm\pmb{\sc Corollary 8.23.} \  {\it Let $n\in\omega$ be 
such that $\omega_n\le\lambda$,
and assume that
$\mu\ge\omega$ if $n=0$.  Then
$\dom_{\omega_n,\lambda^{\aleph_0}}^\mu=\dom_{\omega_n,
\lambda\cdot 2^{\aleph_0}}^\mu$.}\\
\vspace*{0ex}{\it Proof. }The result follows from Proposition 8.22 (0) if
$n\ge 2$, and from Proposition 8.22 (1) if
$n=1$. Let us now turn to the case $n=0$. If $\lambda=\omega$, the 
result is trivial. If $\lambda>\omega$,
we get
\vskip-0,3cm
\centerline{$\dom_{\omega,\lambda^{\aleph_0}}^\mu=\dom_{\omega,2^{\aleph_0}
}^\mu\cdot\dom_{\omega_1,
\lambda}^\mu=\dom_{\omega,2^{\aleph_0}}^\mu\cdot\dom_{\omega,\omega}^\mu
\cdot\dom_{\omega_1,\lambda}^\mu=
\dom_{\omega,2^{\aleph_0}}^\mu \cdot
\dom_{\omega,\lambda}^\mu=\dom_{\omega,\lambda\cdot2^{\aleph_0}}^\mu$}
 
by Propositions 8.22 (1) and 8.15 (0) and Lemma 8.11.\eop

Notice that if $n=0$ and $\mu<\omega$, then 
$\dom_{\omega_n,\lambda^{\aleph_0}}^\mu\not=
\dom_{\omega_n,\lambda\cdot 2^{\aleph_0}}^\mu$ for some values of 
$\lambda$, since
$\dom_{\omega,\lambda^{\aleph_0}}^\mu=
\lambda^{\aleph_0}$ and $\dom_{\omega,\lambda\cdot 
2^{\aleph_0}}^\mu=\lambda\cdot 2^{\aleph_0}$
by Propositions  8.3 and 1.3 (0).\\

\bigskip
\hskip-0,2cm\pmb{\sc Corollary 8.24.} \  {\it If $\lambda\ge 2^{<\kappa}$, then
$\dom_{\kappa,\lambda^{<\kappa}}^\mu=
\dklmu$.}\\
\vspace*{2ex}{\it Proof. }By Proposition 8.22 (1) and Lemma 8.11.\eop

\bigskip
\hskip-0,2cm\pmb{\sc Corollary 8.25.} \  {\it Let $\sigma$ be an 
infinite cardinal such that
$cf(\sigma) < \kappa$ and \break $\kappa^{cf(\sigma)} <\sigma < 
\lambda\leq\sigma^{cf(\sigma)}$. Then \
$\dom_{\kappa,\lambda}^\mu = \dom_{\kappa,\sigma}^\mu$.
}
\\
\vspace*{1ex}{\it Proof. } If $(cf(\sigma))^+<\kappa$, then by Lemma 
8.11 and Proposition 8.22 (0),
\vskip-0,6cm
$$\dom_{\kappa,\lambda}^\mu \leq \dom_{\kappa,\sigma^{cf(\sigma)}}^\mu =
\dom_{\kappa,\sigma\cdot\kappa^{cf(\sigma)}}^\mu = 
\dom_{\kappa,\sigma}^\mu \leq
\dom_{\kappa,\lambda}^\mu.$$
If $(cf(\sigma))^+ = \kappa$, then by Lemma 8.11 and Proposition 8.22 (1),
$$\dom_{\kappa,\lambda}^\mu \leq \dom_{\kappa,\sigma^{cf(\sigma)}}^\mu =
\dom_{\kappa,2^{cf(\sigma)}}^\mu \cdot \dom_{\kappa,\sigma}^\mu = 
\dom_{\kappa,\sigma}^\mu \leq
\dom_{\kappa,\lambda}^\mu . \eqno\hbox{\eop} $$

We finally investigate $\dom_{\kappa,\lambda}^{\mu^{<\rho}}$.\\
 
\hskip-0,2cm\pmb{\sc Proposition 8.26.} \  {\it Let $\rho$ be an 
infinite cardinal with
$\rho\le\mu$.  Then $\dom_{\kappa,\lambda}^{\mu^{<\rho}}$ is the 
least cardinality of any
$F\subseteq{^{(P_\rho(
\mu))}}\!{\pkl}$ with the property that for any 
$g\in{^{(P_\rho(\mu))}}\!{\pkl}$,
there is $f\in F$ with $\{d\in P_\rho(\mu):g(d)\subseteq
f(d)\}\in I_{\rho,\mu}^*$.}\\
\vskip-0,3cm
\vspace*{0ex}{\it Proof. }Let $F\subseteq{^{(P_\rho(\mu))}}\!{\pkl}$
be such that for every $g\in{^{(P_\rho(\mu))}}\!{\pkl}$,
there is $f\in F$ with the property that $\{
d\in P_\rho(\mu):g(d)\subseteq f(d)\}\in I_{\rho,\mu}^*$.
By Corollary 1.5, there are $A_e\in\widehat{e}\cap I_{\rho,\mu}^+$
for $e\in P_\rho(\mu)$  such that (a) $|A_e|=\mu^{<\rho}$ for every
$e\in P_\rho(\mu)$, \break  (b)
$A_e\cap A_{e^\prime}=\emptyset$  whenever $e,e^\prime\in P_\rho(\mu)$ are
  such that
$e\not=e^\prime$, and (c) $\displaystyle\bigcup_{e\in P_\rho(\mu) 
}A_e=P_\rho(\mu)$.
Select  a bijection $j_e:A_e\longrightarrow P_\rho(\mu)$ for each 
$e\in P_\rho(\mu)$.
Given $h\in
{^{(P_\rho(\mu))}}\!{\pkl}$, define $g\in{^{(P_\rho(\mu))}}\!{\pkl}$ so that
$g(d)=h(j_e(d))$ whenever $d\in A_e$. Pick $f\in F$ and $e\in 
P_\rho(\mu)$ so that
$\widehat{e}\subseteq\{d\in P_\rho(\mu):g(d)\subseteq f(d)\}$. Then 
$h(j_e(d))\subseteq
(f\circ j_e^{-1})(j_e(d))$
for every $d\in A_e$. Thus 
$\dom_{\kappa,\lambda}^{\mu^{<\rho}}\le|F|\cdot\mu^{<\rho}$,
and therefore
$\dom_{\kappa,\lambda}^{\mu^{<\rho}}\le|F|$ by Propositions 8.2 and 8.4.\eop

\bigskip
\hskip-0,2cm\pmb{\sc Proposition 8.27.} \  {\it Let $\rho$ be an 
infinite cardinal such
that
$\rho\leq\mu$ and $2^\tau<\kappa$ for every  cardinal $\tau<\rho$. Then
$\dom_{\kappa,\lambda}^{\mu^{<\rho}}=\dom_{\kappa,
\lambda}^{u(\rho,\mu)}$.}
 
\bigskip
\vspace*{0ex}{\it Proof. }We have
$\dom_{\kappa,\lambda}^{u(\rho,\mu)}\le\dom_{\kappa,\lambda}^{\mu^{<\rho}}$ by Lemma 8.11.  For the
other inequality, fix $A\in I_{\rho,\mu}^+$ and $F\subseteq{^A}\!{\pkl}$ with
the property  that for every $g\in {^A}\!{\pkl}$, there is $f\in F$ such that
$g(a)\subseteq f(a)$ for all $a\in A$.  For
$f\in F$, define $f^\prime\in {^{(P_\rho(\mu))}}\!{\pkl}$ as follows: 
given $b\in P_\rho(\mu)$, pick
$a\in A$ with $b\subseteq a$, and set $f^\prime(b)=f(a)$. Now given
$h\in {^{(P_\rho(\mu))}}\!{\pkl}$, define $g\in {^A}\!{\pkl}$ by  $g(a)=
\displaystyle\bigcup_{b\subseteq a} h(b)$. Select $f\in F$ so that 
$g(a)\subseteq f(a)$ for all
$a\in A$. Then
$h(b)\subseteq f^\prime (b)$ for all $b\in P_\rho(\mu)$.\eop

\vskip1,5cm

{\Large\bf 9. $cof(J)$}

\bigskip
This  section is devoted to the computation of $cof(\nsdt)$.

\bigskip
\hskip-0,2cm\pmb{\sc Lemma 9.1.} \  {\it Assume 
$\nabla^\dt\ikl\subseteq J$. Then $cof(
J)\ge\dom_{\kappa,\lambda}^{|P_{\overline{\theta}}(\delta)|}$.}

\medskip
\vspace*{0ex}{\it Proof. }Fix $S\subseteq J$ with
$J=\displaystyle\bigcup_{B\in S}P(B)$. For  $B\in S$, define $h_B:
P_{\overline{\theta}}(\delta)\longrightarrow\pkl-B$ so that
$e\in P_{|\overline{\theta}\cap h_B(e)|}(h_B(e))$ for all
$e\in P_{\overline{\theta}}(\delta)$. Given 
$g:P_{\overline{\theta}}(\delta)\longrightarrow\pkl$, there is $B\in 
S$
with $\pkl-B\subseteq\build{\Delta}_{e\in 
P_{\overline{\theta}}(\delta)}^{}\widehat{g(e)}$
by Proposition 2.2 (0)
and Corollary 2.8 ((iv)$\rightarrow$ (ii)). Then $g(e)\subseteq 
h_B(e)$ for every $e\in
P_{\overline{\theta}}(\delta)$ .\eop

\bigskip
\hskip-0,2cm\pmb{\sc Proposition 9.2.} \  {\it
$cof(\nsdt|A)=\dom_{\kappa,\lambda}^{|P_{\overline{\theta}}(\delta)|}$ 
for every
$A\in(\nsdt)^+$.}\\
\vspace*{0ex}{
\it Proof. }Let us first observe that if
$f:P_{\overline{\theta}\cdot 3}(\delta)\longrightarrow\pkl$ and
$g:P_{\overline{\theta}\cdot 3}(\delta)\longrightarrow\pkl$ are such 
that $f(e)\subseteq
g(e)$ for all $e\in
  P_{\overline{\theta}\cdot 3}(\delta)$, then $C_g^{\kappa,\lambda}\subseteq
C_f^{\kappa,\lambda}$. Hence
$cof(\nsdt)\le\dom_{\kappa,\lambda}^{|P_{\overline{\theta}}(\delta)|}$ 
by Lemma 2.12. So given $A
\in(\nsdt)^+$, we have $cof(\nsdt|A)
\le\dom_{\kappa,\lambda}^{|P_{\overline{\theta}}(\delta)|}$ by 
Proposition 1.6. The reverse inequality holds by Lemma
9.1 since $\nsdt|A$ is \mbox{$\dt$-normal}.\eop

The following is well-known.\\

\bigskip
\hskip-0,2cm\pmb{\sc Corollary 9.3.} \  {\it 
$cof(\ikl|A)=u(\kappa,\lambda)$ for every
$A\in\iklp$.}
\bigskip

\vspace*{0ex}{
\it Proof. }We have $\ikl=NS_{\kappa,\lambda}^{\lbrack 
1\rbrack^{<2}}$ by Propositions 3.5,
2.10 and 2.18 (0).
  So the result follows from Propositions 9.2 and 8.3.\eop

It follows from Proposition 3.5 and Corollary 9.3 that 
$cof(\nsdt|A)=u(\kappa,\lambda)$ for all
$A\in(\nsdt)^+$ if $\delta<\kappa$. For $\delta\ge\kappa$ we have the 
following.\\

\hskip-0,2cm\pmb{\sc Corollary 9.4.} \  {\it Assume $\delta\ge\kappa$. Then
$cof(\nsdt|A)=\dom_{\kappa,\lambda}^{u(\overline{\theta}\cdot\aleph_0, 
|\delta|)}$ for every\break
$A\in(\nsdt)^+$.}
\bigskip

\vspace*{0ex}{
\it Proof. }By Propositions 9.2, 8.26 and 2.18.\eop
\vskip-0,3cm
Under  GCH, we obtain the following values.

\eject
\hskip-0,2cm\pmb{\sc Corollary 9.5.} \  {\it Assume that the} GCH 
{\it holds and $\delta\ge\kappa$,
and let
$A\in(\nsdt)^+$. Then}
\begin{enumerate}
\item[a)] \vskip-0,3cm $cof(\nsdt|A)=\lambda^{++}$ if $\delta=\lambda$ and
$cf(\lambda)<\overline{\theta}$.

\item[b)] \vskip-0,3cm $cof(\nsdt|A)=\lambda^+$ if
$cf(\lambda)\le|\delta|^{<\overline{\theta}}<\lambda$, or
$\lambda\le|\delta|^{<\overline{
\theta}}$ and $cf(\lambda)\ge\overline{\theta}$.

\item[c)] \vskip-0,3cm $cof(\nsdt|A)=\lambda$ if $|\delta|^{<\overline{\theta}}
<cf(\lambda)$.
\end{enumerate}
\vspace*{0ex}{\it Proof. }By Propositions 8.9 and 9.2.\eop

\hskip-0,2cm\pmb{\sc Corollary 9.6.} \  {\it Let $\delta'$ be an ordinal with
$\kappa\leq\delta'\leq\lambda$, and $\theta'$ be a cardinal with 
$2\leq\theta'\leq\kappa$. Assume
that either $\lambda^{|\delta'|^{<\overline{\theta'}}} \leq 
|\delta|^{<\overline\theta}$, or
$\lambda^{|\delta']^{<\overline{\theta'}}} = \lambda$ and $cf(\lambda)\leq
|\delta|^{<\overline\theta}$.  Then there is no $A\in {\big(
NS_{\kappa,\lambda}^{[\delta]^{<\theta}}\big)}^+ \cap {\big(
NS_{\kappa,\lambda}^{[\delta']^{<\theta'}}\big)}^+$  such that \
$NS_{\kappa,\lambda}^{[\delta]^{<\theta}} | A = 
NS_{\kappa,\lambda}^{[\delta']^{<\theta'}} | A$.
}\bigskip

\vspace*{0ex}{\it Proof. } It suffices to observe that by 
Propositions 8.2 and 8.4,  (a) if
$\lambda^{|\delta'|^{<\overline{\theta'}}} \leq 
|\delta|^{<\overline\theta}$,  then
$\dom_{\kappa,\lambda}^{|\delta'|^{<\overline{\theta'}}} \leq 
\lambda^{|\delta'|^{<\overline{\theta'}}}
\leq |\delta|^{<\overline\theta} < 
\dom_{\kappa,\lambda}^{|\delta|^{<\overline\theta}}$, and  (b)
if $\lambda^{|\delta'|^{<\overline{\theta'}}} = \lambda$ and $cf(\lambda) \leq
|\delta|^{<\overline\theta}$, then
$$\dom_{\kappa,\lambda}^{|\delta'|^{<\overline{\theta'}}}
\leq\lambda^{|\delta'|^{<\overline{\theta'}}}= \lambda <
\dom_{\kappa,\lambda}^{|\delta|^{<\overline\theta}}.
\eqno\hbox{\eop}$$

\bigskip
\hskip-0,2cm\pmb{\sc Corollary 9.7.} \  {\it Assume $\delta\geq\kappa$. Then
$${\rm cof}~\big( NS_{\kappa,\lambda}^{[\delta]^{<\theta}}\big) = 
{\rm cof}~ \big( NS_{\kappa,
|\delta|}^{[|\delta|]^{<\theta}}\big) \cdot {\rm cov}\ \big(\lambda,
(|\delta|^{<\overline\theta}\big)^+, \big( 
|\delta|^{<\overline\theta}\big)^+, 2\big).$$}
\bigskip
\vspace*{0ex}{\it Proof. } If $\overline\theta < \kappa$,  then by 
Propositions 8.22 (0) and 2.23,
$\dom_{\kappa, 
|\delta|^{<\overline\theta}}^{|\delta|^{<\overline\theta}} = 
\dom_{\kappa, |\delta|
\cdot \kappa^{<\overline\theta}}^{|\delta|^{<\overline\theta}} = \dom_{\kappa,
|\delta|}^{|\delta|^{<\overline\theta}}.$  \break If \ 
$\overline\theta = \kappa$,  \ then \ by
Propositions 8.22 (1) \ and \ 2.18 (1) \ and \ Lemma 8.11,  \ \ $\dom_{\kappa,
|\delta|^{<\overline\theta}}^{|\delta|^{<\overline\theta}} =$

$\dom_{\kappa, 2^{<\kappa}}^{|\delta|^{<\overline\theta}} \cdot \dom_{\kappa,
|\delta|}^{|\delta|^{<\overline\theta}} = \dom_{\kappa, 
|\delta|}^{|\delta|^{<\overline\theta}}.$
In any case we have $\dom_{\kappa, 
|\delta|^{<\overline\theta}}^{|\delta|^{<\overline\theta}} =
\dom_{\kappa, |\delta|}^{|\delta|^{<\overline\theta}}$. Hence if 
~$|\delta|^{<\overline\theta} <
\lambda$, we can infer from Corollary~8.17 that
$$\dom_{\kappa,\lambda}^{|\delta|^{<\overline\theta}} = d_{\kappa,
|\delta|}^{|\delta|^{<\overline\theta}} \cdot u {\big( \big( 
|\delta|^{<\overline\theta}\big)}^+,
\lambda\big) = \dom_{\kappa, |\delta|}^{|\delta|^{<\overline\theta}} \cdot
\hbox{cov}~\big(\lambda, {(|\delta|^{<\overline\theta})}^+, 
{(|\delta|^{<\overline\theta})}^+,
2\big).$$
If $|\delta|^{<\overline\theta} \geq\lambda$, Lemma 8.11 tells us 
that \  $\dom_{\kappa,
|\delta|}^{|\delta|^{<\overline\theta}} \leq 
\dom_{\kappa,\lambda}^{|\delta|^{<\overline\theta}}
\leq \dom_{\kappa,
|\delta|^{<\overline\theta}}^{|\delta|^{<\overline\theta}}$, \  so
\quad \break$d_{\kappa,\lambda}^{|\delta|^{<\overline\theta}} =
\dom_{\kappa, 
|\delta|}^{|\delta|^{<\overline\theta}} = \dom_{\kappa, 
|\delta|}^{|\delta|^{<\overline\theta}}
\cdot \hbox{cov}~\big(\lambda, {(|\delta|^{<\overline\theta})}^+,
{(|\delta|^{<\overline\theta})}^+, 2
\big)$. \eop

\newpage
\newdimen\margeg \margeg=0pt
\def\bb#1&#2&#3&#4&#5&{\par{\parindent=0pt
     \advance\margeg by 1.1truecm\leftskip=\margeg
     {\everypar{\leftskip=\margeg}\smallbreak\noindent
     \hbox to 0pt{\hss [#1]~~}{#2 - }{\it #3}, {#4.}\par\medskip
     #5 }
\medskip}}
{\bf\Large References}
\bigskip

\bb A&Y. ABE&A hierarchy of filters smaller than $CF_{\kappa\lambda}$&
          Archive for Mathematical Logic 36 (1997), 385-397& &

\bb C&D.M. CARR&The minimal normal filter on $P_\kappa \lambda$&
         Proceedings of the American Mathematical Society 86 (1982),
         316 - 320& &

\bb CLP&D.M. CARR, J.P. LEVINSKI and D.H. PELLETIER&On the existence 
of strongly normal
         ideals on $P_\kappa \lambda$&Archive for Mathematical Logic 30
        (1990), 59 - 72& &

\bb CP&D.M. CARR and D.H. PELLETIER&Towards a structure theory for ideals
        on $\pkl$&in~:   Set Theory and its Applications (J. 
Stepr\={a}ns and S. Watson, eds.),
Lecture Notes in  Mathematics 1401, Springer, Berlin, 1989, pp. 41 - 54& &


\bb DoM&H.D. DONDER and P. MATET&Two cardinal versions of 
diamond&Israel Journal of Mathematics 83  (1993), 1 - 43& &

\bb D&M. D\v ZAMONJA&On $P_\kappa\lambda$-combinatorics using a third cardinal&
          Radovi Matemati\v cki 9 (1999), 145-155& &

\bb EHM\'aR&P. ERD\"OS, A. HAJNAL, A. M\'AT\'E and R. 
RADO&Combinatorial Set Theory~:
          Partition Relations for Cardinals&Studies in Logic and the 
Foundations of Mathematics
vol. 106, North-Holland, Amsterdam, 1984& &

\bb F&Q. FENG&An ideal characterization of Mahlo cardinals&Journal of 
Symbolic Logic 54 (1989), 467 - 473& &

\bb HoStW&M. HOLZ, K. STEFFENS and E. WEITZ&Introduction to Cardinal 
Arithmetic&Birkh\"auser, Basel,
1999& &

\bb J&T.J. JECH&The closed unbounded filter over 
$P_\kappa(\lambda)$&Notices of the American Mathematical Society 18
(1971), 663& &


\bb K&K. KUNEN&Set Theory&North-Holland, Amsterdam, 1980& &

\bb L&A. LANDVER&Singular Baire numbers and related topics&Ph. D.
Thesis, University  of Wisconsin-Madison, 1990& &


\bb M1&P. MATET&Un principe combinatoire en relation avec
l'ultranormalit\'e des id\'eaux&Comptes Rendus de l'Acad\'emie des Sciences de
Paris, S\'erie I, 307 (1988), 61 - 62& &

\bb M2&P. MATET&Concerning stationary subsets of $\lk$& in~: Set Theory
and its Applications (J. Stepr\={a}ns and S. Watson, eds.), Lecture 
Notes in Mathematics 1401,
Springer, Berlin, 1989, pp. 119 - 127& &

\bb M3&P. MATET&Partition relations for $\kappa$-normal ideals on
$P_\kappa(\lambda)$&Annals of Pure and Applied Logic 121 (2003), 89-111& &

\bb M4&P. MATET&Covering for category and combinatorics on \ 
$P_\kappa(\lambda )$
&preprint& &

\bb MP\'eS&P. MATET, C. P\'EAN and S. SHELAH&Cofinality of normal 
ideals on $P_\kappa
(\lambda)$ II&Israel Journal of Mathematics, to appear& &

\bb Me&T.K. MENAS&On strong compactness and supercompactness&Annals 
of Mathematical Logic 7
(1974),  327 - 359& &

\bb S1&S. SHELAH&Cardinal Arithmetic&Oxford Logic Guides vol. 29, 
Oxford University Press,
Oxford, 1994& &

\bb S2&S. SHELAH&The generalized continuum hypothesis revisited&Israel
Journal of Mathematics 116 (2000), 285-321& &

\vfill\eject

\end{document}